%% file: Brooks-beyond.tex
\documentclass[12pt]{amsart}
\usepackage{amsmath, amsthm, amssymb}
\usepackage[top=1.25in, bottom=1.25in, left=1.0in, right=1.0in]{geometry}
\usepackage{hyperref}
\usepackage{mathscinet}
\usepackage{tikz,tkz-graph,tkz-berge}
\usetikzlibrary{positioning}
\usepackage{verbatim, xcolor}

\usetikzlibrary{arrows}
\usetikzlibrary{decorations.markings}
\newcommand{\boundellipse}[3]
{(#1) ellipse (#2 and #3)
}

\makeatletter
\newtheorem*{rep@theorem}{\rep@title}
\newcommand{\newreptheorem}[2]{
\newenvironment{rep#1}[1]{
 \def\rep@title{#2 \ref{##1}}
 \begin{rep@theorem}}
 {\end{rep@theorem}}}
\makeatother

\theoremstyle{plain}
\newtheorem{thm}{Theorem}
\newreptheorem{thm}{Theorem}
\newtheorem*{Brooks}{Brooks' Theorem}
\newtheorem*{KernelLemma}{Kernel Lemma}
\newtheorem*{BK}{Borodin--Kostochka Conjecture}
\newtheorem*{Reed}{Reed's Conjecture}
\newtheorem*{AT-thm}{Alon--Tarsi Theorem}
\newtheorem*{Rubin}{Rubin's Block Lemma}

\newtheorem{lem}[thm]{Lemma}
\newreptheorem{lem}{Lemma}

\newreptheorem{cor}{Corollary}
\theoremstyle{definition}

\newcommand{\IN}{\mathbb{N}}

\newcommand{\LB}{{L}_B}

\newcommand{\chil}{{\chi_{\ell}}}
\newcommand{\chiol}{{\chi_{OL}}}
\newcommand{\DeltaG}{{\Delta}}
\newcommand{\erdos}{Erd\H{o}s}
\newcommand{\CondOne}{\hyperref[cond1]{Condition (1)}}
\newcommand{\CondTwo}{\hyperref[cond2]{Condition (2)}}

\newcommand{\set}[1]{\left\{ #1 \right\}}

\newcommand{\card}[1]{\left|#1\right|}
\newcommand{\size}[1]{\left\Vert#1\right\Vert}
\newcommand{\ceil}[1]{\left\lceil#1\right\rceil}

\newcommand{\func}[3]{#1\colon #2 \rightarrow #3}

\newcommand{\irange}[1]{\{1,\ldots,#1\}}

\newcommand{\brackets}[1]{\left[ #1 \right]}
\newcommand{\DefinedAs}{=}

\begin{document}
\title{Brooks' Theorem and Beyond}
\author{Daniel W. Cranston}\thanks{Virginia Commonwealth University, Richmond,
VA.  \texttt{dcranston@vcu.edu}} 
\author{Landon Rabern}\thanks{LBD Data Solutions, Orange, CT. 
\texttt{landon.rabern@gmail.com}}
\begin{abstract}
We collect some of our favorite proofs of Brooks' Theorem, highlighting
advantages and extensions of each.  
The proofs illustrate some of the major techniques in graph coloring, such as
greedy coloring, Kempe chains, hitting sets, and the Kernel Lemma. We also
discuss standard strengthenings of vertex coloring, such as list coloring,
online list coloring, and Alon--Tarsi orientations, since analogues of Brooks'
Theorem hold in each context.  We conclude with two conjectures along the lines
of Brooks' Theorem that are much stronger, the Borodin--Kostochka Conjecture
and Reed's Conjecture.
\end{abstract}
\date{3 March 2014}
\maketitle


\bigskip

Brooks' Theorem is among the most fundamental results in graph coloring.
In short, it characterizes the (very few) connected graphs for which an obvious
upper bound on the chromatic number holds with equality. 
It has been proved and reproved using a wide range of techniques, and the
different proofs  generalize and extend in many directions. 
In this paper we share some of our favorite
proofs. 
In addition to surveying Brooks' Theorem, we aim to illustrate many of the
standard techniques in vertex coloring\footnote{One common coloring approach
notably absent from this paper is the probabilistic method.  Although the
so-called naive coloring procedure has been remarkably effective in proving
numerous coloring conjectures asymptotically, we are not aware of any
probabilistic proof of Brooks' Theorem.  For the reader interested in this
technique, we highly recommend the monograph of Molloy and Reed ``Graph
Colouring and the Probabilistic Method'' \cite{MolloyR-GCPM}.};  furthermore,
we prove versions of Brooks' Theorem for standard strengthenings of vertex
coloring, including list coloring, online list coloring, and Alon--Tarsi
orientations.  We present the proofs roughly in order of increasing complexity,
but each section is self-contained and the proofs can be read in any order.
Before we state the theorem, we need a little background.

A \emph{proper coloring} assigns colors, denoted by positive
integers, to the vertices of a graph so that endpoints of each edge get
different colors.  A graph $G$ is \emph{$k$-colorable} if it has a proper
coloring with at most $k$ colors, and its \emph{chromatic number} $\chi(G)$ is
the minimum value $k$ such that $G$ is $k$-colorable.  If a graph $G$ has
maximum degree $\Delta$, then $\chi(G)\le \Delta+1$, since we can repeatedly
color an uncolored vertex with the smallest color not already used on its
neighbors.  Since the proof of this upper bound is so easy, it is natural to
ask whether we can strengthen it.  The answer is yes, nearly always.

A \emph{clique} is a subset of vertices that are pairwise adjacent.
If $G$ contains a clique $K_k$ on $k$ vertices, then 
$\chi(G)\ge k$, since all clique vertices need distinct colors.
Similarly, if $G$ contains an odd length cycle, then $\chi(G)\ge 3$, even when
$\Delta=2$.
In 1941, Brooks~\cite{Brooks41} proved that these are the only two cases in
which we \emph{cannot} strengthen our trivial upper bound on $\chi(G)$.
 
\begin{Brooks}
Every graph $G$ with maximum
degree $\Delta$ has a $\Delta$-coloring unless either (i) $G$ contains
$K_{\Delta+1}$ or (ii) $\Delta=2$ and $G$ contains an odd cycle.  
\end{Brooks}
%
We typically emphasize the case $\Delta\ge 3$, and often write: ``Every
graph $G$ satisfies $\chi \le \max\{3,\omega,\Delta\}$,''
where $\omega$ is the size of the largest clique in $G$.
This formulation follows Brooks' original paper, and it has the
benefit of 
saving us from considering odd cycles in every proof.

Before the proofs, we introduce our notation, which is fairly standard.  The
\emph{neighborhood} $N(v)$ of a vertex $v$ is the set of vertices adjacent to
$v$.  The \emph{degree} $d(v)$ of a vertex $v$ is $\card{N(v)}$.
For a graph $G$ with vertex set $V(G)$, we often write $\card{G}$ to denote
$\card{V(G)}$.
The \emph{clique number} $\omega(G)$ of $G$ is the size of its largest clique.
We denote the maximum and minimum degrees of $G$ by $\Delta(G)$ and $\delta(G)$.
We write simply $\omega$, $\Delta$, $\delta$, or $\chi$ when refering
to the original graph $G$, rather than to any subgraph.
The subgraph of $G$ induced by vertex set $V_1$ is denoted $G[V_1]$, and 
the subgraph induced by a clique is \emph{complete}. 
To \emph{color greedily} is to consider the vertices
in some order, and color each vertex $v$ with the smallest color
not yet used on any of its neighbors.   For every graph, there exists a vertex
order under which greedy coloring is optimal (given an optimal coloring,
consider the vertices in order of increasing color).  However,
considering each of the $n!$ vertex orderings for an $n$-vertex graph is
typically impractical.


In our proofs, we often assume that a counterexample exists, and this assumption
leads us to a contradiction.
A \emph{minimum counterexample} $G$ to Brooks Theorem is one minimizing $|G|$.
To avoid repetition later, we note the following here, which is useful
in multiple proofs.  Any minimum counterexample $G$ must satisfy $\chi >
\max\{\omega, \Delta\}$, and thus $\chi = \Delta + 1$.  
If $H$ is a proper induced subgraph of $G$, then minimality of $G$ gives
$\chi(H) \le \max\{3,\omega(H),\Delta(H)\} \le \Delta$.
So for every vertex $v$, $G-v$ is $\Delta$-colorable; since $G$ is not
$\Delta$-colorable, $v$ must have a neighbor in all $\Delta$ color classes in
any $\Delta$-coloring of $G-v$.  Since this is true for every vertex $v$, it
follows that $G$ must be $\Delta$-regular.


\section{Greedy coloring}
\label{lovasz}
To motivate our first two proofs, we return to the observation that every graph
satisfies $\chi\le \Delta+1$.  The idea is to color the vertices greedily in
an arbitrary order.  To prove $\chi \le \Delta$, it suffices to order the
vertices so that at the point when each vertex gets colored, it still has an
uncolored neighbor.  Of course, this is impossible, since in any order
the final vertex will
have no uncolored neighbors.  Nonetheless, we \emph{can} find an order such
that every vertex except the last still has an uncolored neighbor when it gets
colored.  Such an order yields a $\Delta$-coloring of $G-v$ for each choice
of a final vertex $v$.  Our first two proofs show two different ways to ensure
that we can extend this $\Delta$-coloring to $G$.
\bigskip

In 1975, Lov\'{a}sz published a 3-page article 
\cite{Lovasz1975269} 
entitled ``Three Short Proofs in Graph Theory.''
It included the following proof of Brooks' Theorem 
by coloring greedily in a good order. 
The proof needs a few notions of connectedness. 
A \emph{cutset} is a subset $V_1\subset V$ such that $G\setminus V_1$ is
disconnected.  If a single vertex is a cutset, then it is a
\emph{cutvertex}.  A graph is \emph{$k$-connected} if every cutset has size at
least $k$.
A \emph{block} is a maximal 2-connected subgraph.  The \emph{block graph} of a
graph $G$ has the blocks of $G$ as its vertices and has two blocks adjacent if
they intersect.  It is easy to see that every block graph is a forest; each
leaf of a block graph corresponds to an \emph{endblock} in the original graph.

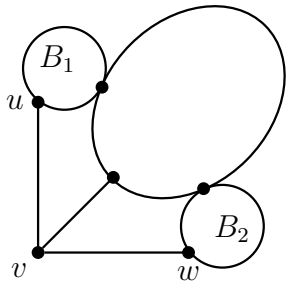
\begin{figure}[!ht]
\begin{center}

\begin{tikzpicture}[line width=.03cm]
\tikzstyle{vert}=[draw,shape=circle,fill=black,minimum size=4pt, inner sep=0pt] 

\draw (0,0) node[vert,label={[label distance=-.1cm]225:$v$}] (v) {} to
      (0,2) node[vert,label={[label distance=-.05cm]180:$u$}] (u) {};
\draw (v) to  
      (2,0) node[vert,label={[label distance=-.05cm]270:$w$}] (w) {};

\draw \boundellipse{.35,2.45}{.55}{.55}; 
\draw \boundellipse{2.45,.35}{.55}{.55}; 
\draw[rotate around={45:(2.,2.)}] (2.,2.) ellipse (41pt and 31pt);

\draw (2.2,0.85) node[vert] {};
\draw (0.85,2.2) node[vert] {};
\draw (v) to (1.0,1.0) node[vert] {};

\draw (.8,2.58) node[label=left:{$B_1$}] {};
\draw (2.58,.8) node[label=below:{$B_2$}] {};
\end{tikzpicture}

\end{center}
\caption{$G$ contains an induced path $uvw$ such that $G\setminus\{u,w\}$ is
connected.}
\label{fig1}
\end{figure}

\begin{lem}\label{TwoConnectedHasGoodP3}
Let $G$ be a $2$-connected graph with $\delta(G) \geq 3$.
If $G$ is not complete, then $G$ contains an induced path on 3 vertices, say
$uvw$, such that $G \setminus\{u, w\}$ is connected.
\end{lem}
\begin{proof}
Since $G$ is connected and not complete, it contains an induced path on 3
vertices. If $G$ is $3$-connected, any such path will do.  Otherwise,
let $\set{v, x} \subset V(G)$ be a cutset.  Since $G - v$ is not
$2$-connected, it has at least two endblocks $B_1, B_2$.  Since $G$ is
$2$-connected, each endblock of $G-v$ has a noncutvertex adjacent to $v$ (see
Figure~\ref{fig1}).  Let $u \in B_1$ and $w \in B_2$ be such vertices.  Now $G
\setminus\{u, w\}$ is connected since $d(v) \geq 3$.  So $uvw$ is our desired
induced path.
\end{proof}

\begin{proof}[Proof~\ref{lovasz}  of Brooks' Theorem.]
Let $G$ be a connected graph.
First suppose that $G$ has a vertex $v$ with $d(v)<\Delta$.  
We color greedily in order of decreasing distance to $v$ (breaking ties
arbitrarily).  For each vertex $u$ other than $v$, when $u$ gets colored, some
neighbor $w$ on a shortest path in $G$ from $u$ to $v$ is uncolored, so we
use at most $\Delta$ colors.  Since $d(v)<\Delta$, we can color $v$ last.  
Similarly, suppose $G$ has a cutvertex $v$. 
Now for each component $H$ of $G-v$, we can $\Delta$-color $H+v$, since $v$ has
fewer than $\Delta$ neighbors in $H$.  By permuting the color classes to agree
on $v$, we get a $\Delta$-coloring of $G$.  

Now assume that $G$ is $\Delta$-regular, $2$-connected, and not complete. 
Let $uvw$ be the induced path guaranteed by Lemma
\ref{TwoConnectedHasGoodP3}.
Color $u$ and $w$ with color 1; now as before, color the remaining vertices
greedily in order of decreasing distance in $G\setminus\{u, w\}$ from $v$.  
Again we can color $v$ last, this time because it has
two neighbors with the same color.
\end{proof}

Lov\'{a}sz's proof has many variations.  
Bondy \cite{bondy2003short} used a depth-first-search tree to construct the
path in Lemma \ref{TwoConnectedHasGoodP3},
and Bryant \cite{bryant1996characterisation} used yet another method.
Schrijver's proof 
\cite{schrijver2003combinatorial} 
skips Lemma \ref{TwoConnectedHasGoodP3} by
using greedy coloring only for $3$-connected graphs and handling two-vertex
cutsets by patching together colorings of the components.

\section{Kempe chains}
\label{kempe}
The most famous theorem in graph theory is the 4 Color Theorem: Every
planar graph is 4-colorable.
In 1852, Guthrie asked whether this was true.   Two years
later the problem appeared in \emph{Athen\ae um}~\cite{Guthrie1854, McKay14},
a London literary journal, where it
attracted the attention of mathematicians.
In 1879, Kempe published a proof.
Not until 1890 did Heawood highlight a flaw in Kempe's purported
solution.\footnote{In 1880, Tait published a second proof.  But, alas, it, too,
was founding wanting.}  Fortunately, Heawood largely salvaged Kempe's ideas,
and proved the 5~Color Theorem.  The key tool in this work is now called a
Kempe chain.  This technique is among the most well-known in graph coloring. 
It yields a short proof that every bipartite graph has a proper
$\Delta$-edge-coloring, and similar ideas show that \emph{every} graph has a
proper $(\Delta+1)$-edge-coloring (see Section 5.3 of Diestel
\cite{DiestelBook}). 

In 1969 Mel'nikov and Vizing \cite{MelnikovV69} used Kempe chains to give the
following elegant proof of Brooks' Theorem.
We phrase this proof in terms of a minimal counterexample $G$, and for
an arbitrary vertex $v$, we color $G-v$ by minimality.  To turn the proof into
an algorithm, we can simply color greedily toward $v$, as in our first proof;
thus, the ``hard part'' is once again showing how to color this final vertex
$v$.

\begin{center}
\begin{figure}[h!tb]
\begin{tikzpicture}[line width=.03cm,scale=.8]
\tikzstyle{vert}=[draw,shape=circle,fill=black,minimum size=4pt, inner sep=0pt] 
\tikzstyle{small}=[draw,shape=circle,fill=black,minimum size=2pt, inner sep=0pt] 
\begin{scope}[rotate=20]
\foreach \deg in {120,165,210,360}  
{
\draw ({3*cos \deg}, {2*sin \deg}) node[small] (v\deg) {}; 
}

\draw (3*cos 90,2*sin 90) 
-- (v120) node [vert,label={[label distance=-.1cm]150:$v_i$}] {}
-- (v165) node [vert,label={[label distance=-.1cm]180:$v$}] {}
-- (v210) node [vert,label={[label distance=-.05cm]280:$v_j$}] {}
-- (3*cos 245,2*sin 245);
\draw (3*cos 25,2*sin 25) 
-- (v360) node [vert,label={[label distance=-.1cm]280:$u$}] {}
-- (3*cos 320,2*sin 320)
(v360) -- ++(cos 0,sin 0);
\draw (-.4,.2) node (Cij) {$C_{i,j}$};
\draw (1.4,1.6) node (cdots) {\LARGE{$\cdots$}};
\draw (-1.6,-1.4) node (myAnchor) {};
\node [right = of myAnchor, rotate=28] (cdots) {\LARGE{$\cdots$}};

\end{scope}
\end{tikzpicture}
\hspace{.6in}
%
%
\begin{tikzpicture}[line width=.03cm,scale=.6]
\tikzstyle{vert}=[draw,shape=circle,fill=black,minimum size=4pt, inner sep=0pt] 
\tikzstyle{small}=[draw,shape=circle,fill=black,minimum size=2pt, inner sep=0pt] 

\foreach \deg in {180,360}  
{
\draw ({3*cos \deg}, {.9*sin \deg}) node[small] (v\deg) {}; 
}
\draw (0,0.7) node (cdots1) {\LARGE{$\cdots$}};
\draw (0,-.7) node (cdots1) {\LARGE{$\cdots$}};

\draw (3*cos 125,.75*sin 125)  
-- (v180) node[vert,label={[label distance=-.12cm]220:$v_i$}] (vi) {} 
-- (3*cos 235,.75*sin 235);

\draw (3*cos 305,.75*sin 305)  
-- (v360) node [vert,label={[label distance=-.1cm]0:$u$}] {} 
-- (3*cos 55,.75*sin 55);

\draw (.1,2.8) node (cdots1) {\LARGE{$\cdots$}};

\draw (.1,-2.8) node (cdots1) {\LARGE{$\cdots$}};

\draw (-2.0,2.75) -- 
(-3,2.5) node[vert] {} 
-- (-3,2.5) node[vert,label={[label distance=-.08cm]140:$v_j$}] (vj) {} 
-- (-4.5,0) node[vert,label={[label distance=-.08cm]180:$v$}] (v) {} 
-- (-3,-2.5) node[vert,label={[label distance=-.12cm]220:$v_k$}] (vk) {} 
-- (-2.0,-2.75);
\draw (v) -- (vi);

\draw (-1.4,1.6) node (Cik) {$C_{i,j}$}; 
\draw (-1.7,-1.6) node (Cjk) {$C_{i,k}$}; 

\draw (2.15,2.4) -- (3,0) -- (2.15,-2.4);

\end{tikzpicture}

\caption{The left figure shows Claim 2 and the right figure shows
Claim 3 in Proof \ref{kempe} of Brooks' Theorem.}
\end{figure}
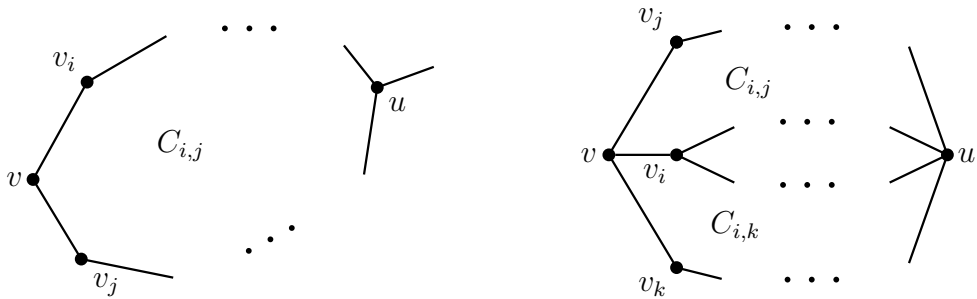
\end{center}

For a proper coloring of a graph $G$, an \emph{$(i,j)$-Kempe chain} is a 
component of the subgraph of $G$ induced by the vertices of colors $i$ and
$j$. 
A \emph{swap} in an $(i,j)$-Kempe chain $H$ swaps the colors on $H$; 
each vertex in $H$ colored $i$ is recolored $j$ and vice versa.  Such a swap
yields another proper coloring.  


\begin{proof}[Proof~\ref{kempe}  of Brooks' Theorem.]
Suppose the theorem is false and let $G$ be a minimum counterexample.
Choose an arbitrary vertex $v$.
By minimality, $G$ is $\Delta$-regular.  Further, in each $\Delta$-coloring
of $G-v$, each color appears on some neighbor of $v$.  Fix an arbitary
$\Delta$-coloring $C$ of $G-v$.  For each
$i\in\irange{\Delta}$, let $v_i$ be the neighbor of $v$ using color $i$.
By a similar argument, for each $v_i$, each color other than $i$ appears on a
neighbor of $v_i$; for otherwise we could recolor $v_i$ and color $v$ with $i$.
For each pair of colors $i$ and $j$, let $C_{i,j}$ denote the $(i,j)$-Kempe
chain containing $v_i$.

\textbf{Claim 1.} \textit{For all $i$ and $j$, $C_{i,j}=C_{j,i}$.}
If not, then after a swap in $C_{i,j}$, two neighbors of $v$ use $j$ and none
use $i$, so we can color $v$ with $i$.

\textbf{Claim 2.} \textit{For all $i$ and $j$, the $(i,j)$-Kempe chain
$C_{i,j}$ containing $v_i$ and $v_j$ is a path.}
If not, then $C_{i,j}$ has a vertex of degree at least 3 (since $v_i$
and $v_j$ have degree 1); let $u$ be the unique such vertex in $C_{i,j}$ that is
closest to $v_i$.  At most $\Delta-2$ colors appear on neighbors of
$u$, so we can recolor $u$.  Now $i$ and $j$ violate Claim~1.

\textbf{Claim 3.} \textit{For all $i$, $j$, $k$, we have $C_{i,j}\cap
C_{i,k}=v_i$.}
Suppose not and choose $u\in C_{i,j}\cap C_{i,k}$, with $u\ne v_i$.  Since $u$
has color $i$, colors $j$ and $k$ each appear on two neighbors of $u$;
so we can recolor $u$.  Now $i$ and $j$ (and also $i$ and $k$) violate
Claim~1.

\textbf{Claim 4.} \textit{Brooks' Theorem is true.}
If the neighbors of $v$ form a clique, then $G=K_{\Delta+1}$ and there is
nothing to prove.  So instead there exist some nonadjacent $v_i$ and $v_j$, say
$v_1$ and $v_2$ by symmetry.  Let $u$ be the neighbor of $v_1$ in $C_{1,2}$.
Now perform a swap in $C_{1,3}$. 
Call this new coloring $C'$, and define $v'_i$ and $C'_{i,j}$ analogously to
$v_i$ and $C_{i,j}$ for $C$.
Since $u$ still uses color 2, clearly $u \in C'_{2,3}$.  By Claim 3, the swap on
$C_{1,3}$ did not disrupt $C_{1,2}$ except at $v_1$; so $u\in C'_{1,2}$.
Now $u\in C'_{1,2}\cap C'_{2,3}$, which violates Claim 3, and gives the desired
contradiction.
\end{proof}

Kostochka and Nakprasit \cite{kostochka2005equitable} took this Kempe chain
proof further by showing that when extending the coloring to $v$, 
we can ensure that only one color class changes size.\footnote{
They proved the following.
For $k \ge 3$, let $G$ be a $K_{k+1}$-free graph with $\Delta(G) \le k$.  Suppose that $G-v$ has a $k$-coloring with color classes $M_1, \ldots, M_k$.  
Then $G$ has a $k$-coloring with color classes $M_1', \ldots, M_k'$ such that $|M_i| \ne |M_i'|$ for exactly one $i$.
\label{KostochkaN}
}
%
Catlin \cite{catlin1979brooks} proved that the $\Delta$-coloring given by
Brooks' Theorem can be chosen so that one of the color classes is a maximum
independent set, and this result is a quick corollary of
the theorem of Kostochka and Nakprasit.
Much earlier, Mitchem \cite{mitchem1978short} gave a
short proof of Catlin's result by modifying an existing $\Delta$-coloring via
Kempe chains.  

\section{Reducing to the cubic case}
\label{cubic}
A natural idea for proving Brooks' Theorem is induction on the maximum degree.
For a graph $G$ with maximum degree $\Delta$ and clique number at most
$\Delta$, suppose we have some independent set $I$ such that $G-I$ has
maximum degree and clique number each at most $\Delta-1$.  If we can color
$G-I$ with $\Delta-1$ colors, then we can extend the coloring to $G$ with one
extra color.  This approach forms the basis for our proofs in the next two
sections.  Of course we must provide a base case for the induction, and we must
also show how to find this very useful set $I$.
\bigskip

Before again proving Brooks' Theorem, we give an easy lemma, which covers
the base case in our inductive proof.
We often prove coloring results by repeatedly extending a partial coloring.
During this process, the lists of valid colors for two uncolored vertices may
differ, depending on the colors already used on their neighbors.  This
motivates the notion of \emph{list coloring}.   Later we develop this idea
further, but for now we need only the following lemma.

\begin{lem}
\label{cycle-lemma}
If each vertex of a cycle $C$ has a list of 2 colors, then $C$ has a proper
coloring from its lists unless $C$ has odd length and all lists are identical. 
\end{lem}
\begin{proof}
Denote the vertices by $v_1,\ldots,v_n$ and let the lists be as
specified.  If all lists are identical and $n$ is even, then we 
alternate colors on $C$.  So suppose that two lists differ.  This implies that
the lists differ on some pair of adjacent vertices, say on $v_1$ and $v_n$. 
Color $v_1$ with a color not in the list for $v_n$.  Now color the vertices
greedily in order of increasing index.
\end{proof}

The following proof is due to the second author~\cite{rabern2012different}.  In
some ways it is simpler than the first two proofs, since this one needs neither
connectivity concepts nor recoloring arguments.

\begin{proof}[Proof~\ref{cubic}  of Brooks' Theorem.]
Suppose the theorem is false and let $G$ be a minimum counterexample.
Recall that $G$ must be $\Delta$-regular.

First, suppose $G$ is $3$-regular.  
A \emph{diamond} is $K_4$ minus an edge.
If $G$ contains an induced diamond $D$, then by minimality we $3$-color $G-D$. 
The two nonadjacent vertices in $D$ each still
have two colors available, so we color them with a common color, and then
finish the coloring.
So $G$ cannot contain diamonds. Since $\delta(G)\ge 2$, $G$ contains an
induced cycle $C$.  Each vertex of $C$ has one neighbor outside of $C$.
Since $\Delta=3$ and $G$ does not contain $K_4$, two vertices of $C$ have
distinct neighbors outside of $C$; call the neighbors $x$ and $y$ (see
Figure~\ref{cubic-fig}).  When $x$
and $y$ are adjacent, let $H \DefinedAs G - C$; otherwise, let $H \DefinedAs
(G-C) + xy$.  Since $G$ does not contain diamonds, $H$ does not contain $K_4$. 
Since $G$ is minimum, $H$ is $3$-colorable.  That is, $G - C$ has a 3-coloring
where $x$ and $y$ get different colors.  
Each vertex of $C$ loses one color to its neighbor outside of $C$, and so still
has two colors available.  Since $x$ and $y$ use different colors, 
by Lemma \ref{cycle-lemma} we
can extend the coloring to $V(C)$, and hence to all of $G$.

Now instead suppose $\Delta \ge 4$.  Since $G$ is not complete, 
it has an induced 3-vertex path, $uvw$. 
By minimality, $G-v$ has a $\Delta$-coloring; 
choose a color class $I$ of this $\Delta$-coloring with $u,w \not \in I$.
Now $\omega(G-I) \le \Delta - 1$; 
this is because any $K_{\Delta}$ in $G-I$ would have to
contain $v$ and all of its neighbors in $G-I$, but its neighbors $u$ and $w$ are
nonadjacent.
Form $I'$ by expanding $I$ to a maximal independent
set, and let $H= G-I'$.  
Since $I'$ is maximal, each vertex in $H$ has a neighbor in $I'$, so
$\Delta(H)\le \Delta-1$.  
If $\Delta(H) = \Delta-1$, then $\omega(H) \le \Delta(H)$, so 
Brooks' Theorem holds for $H$, and $\chi(H)\le \Delta-1$.
Otherwise, $\Delta(H) \le \Delta-2$.  
Now a greedy coloring gives $\chi(H) \le \Delta(H)+1 \le 
\Delta-1$.  
In each case $\chi(H)\le \Delta-1$; now we use
one more color on $I'$ to get $\chi(G)\le \chi(H)+1\le \Delta$.
\end{proof}

\begin{center}
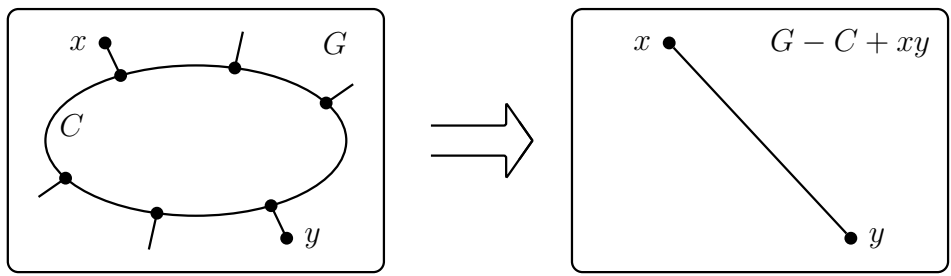
\begin{figure}[htb]

\begin{tikzpicture}[line width=.03cm,scale=.5]
\tikzstyle{vert}=[draw,shape=circle,fill=black,minimum size=4pt, inner sep=0pt] 
\draw \boundellipse{0,0}{4}{2}; 
\foreach \deg in {30,75,120,210,255,300}  
{
\draw ({4*cos \deg}, {2*sin \deg}) node[vert] {} -- ++
({.83*cos \deg}, {sin \deg}); 
}

\foreach \deg/\name\dir in {120/x/left,300/y/right} 
\draw ({4.83*cos \deg}, {3*sin \deg})  node[vert,label=\dir:{$\name$}] {};

\draw[rounded corners] (5,0) node (a) {} -- (5,3.5) -- (-5,3.5) -- (-5,-3.5) -- (5,-3.5)
-- cycle;

\draw (3.7,2.6) node (G) {$G$};
\draw (-3.3,.4) node (C) {$C$};

\begin{scope}[xshift=.25cm] 
\draw[rounded corners=.3mm] (6,.4) -- (8,.4) -- (8,1.0) -- (8.7,0);
\draw[rounded corners=.08mm] (8.35,.5) -- (8.7,0) -- (8.35,-.5);
\draw[rounded corners=.3mm] (8.7,0) -- (8,-1.0) -- (8,-.4) -- (6,-.4);
\end{scope}

\begin{scope}[xshift=15cm]  
\foreach \deg/\name\dir in {120/x/left,300/y/right}
\draw ({4.83*cos \deg}, {3*sin \deg}) node[vert,label=\dir:{$\name$}] (\name){};

\draw (x) -- (y);  

\draw[rounded corners] (5,3.5) -- (-5,3.5) -- (-5,0) node (b) {} -- (-5,-3.5) -- (5,-3.5) -- cycle;

\draw (2.4,2.6) node (G) {$G-C+xy$};

\end{scope}
\end{tikzpicture}

\caption{How to $\Delta$-color a graph with no $K_{\Delta+1}$ when $\Delta=3$.}
\label{cubic-fig}
\end{figure}
\end{center}

In one variation of Proof~\ref{cubic} we do not reduce to the cubic
case \cite{yetanother}.  Instead, we note that, similarly to Lemma
\ref{cycle-lemma}, a $K_k$ has a proper coloring from lists of size $k-1$
unless all the lists are identical.  Hence if $G$ contains a $K_{\Delta}$ (or
an odd cycle when $\Delta = 3$), then we get a $\Delta$-coloring of $G$ by
minimality in the same way as in the $3$-regular case of Proof~\ref{cubic}. 
Otherwise, removing any maximal independent set yields a smaller
counterexample.  

A \emph{hitting set} is an independent set that intersects every maximum
clique.  The reduction to the cubic case in the previous proof is an immediate
consequence of more general lemmas on the existence of hitting sets
\cite{kostochkaRussian,
rabernhitting, KingHitting, tverberg1983brooks}.
Schmerl~\cite{schmerl1982effective} 
extended Brooks' Theorem to all locally finite graphs,
by constructing a recursive hitting set\footnote{By \emph{recursive
hitting set}, we mean an infinite set $S$ of vertices, which is both a hitting
set and a \emph{recursive} set; more formally, $S$ is recursive if for any
vertex $v$, we can decide whether $v$ is in $S$ in finite time.}. 
%
Tverberg also modified his
earlier proof \cite{tverberg1983brooks} to give a shorter constructive proof
\cite{tverberg1984schmerl} of Brooks' Theorem for locally finite graphs.  We
will see this earlier proof in Section~\ref{ktrees}.

\section{k-trees}
\label{ktrees}
Tverberg used $k$-trees \cite{tverberg1983brooks} to give a short proof of
Brooks' Theorem.  Similar to the proof in Section~\ref{cubic}, this
one inductively colors $G-I$ by minimality, where $I$ is a hitting set.  The
differences in the two proofs lie in how we find $I$ and how we handle the base
case.

We define \emph{$k$-trees} as follows.  For $k = 3$,
an odd cycle is a $k$-tree and for $k \ge 4$, a $K_k$ is a $k$-tree. 
Additionally, any graph formed by adding an edge between vertices of degree 
$k-1$ in disjoint $k$-trees is again a $k$-tree.  
For convenience, we write $T_k$ to mean $K_k$ when $k\ge 4$ and to mean an odd
cycle when $k=3$.

The name $k$-tree comes from the fact that contracting each $T_k$ to a single
vertex yields a tree $T$.  A \emph{leaf} in a $k$-tree is a $T_k$ corresponding
to a leaf in $T$.  Each leaf in a $k$-tree has only one vertex of degree $k$ in
$G$, so it is easy to show by induction that each $k$-tree other than $T_k$ has
at least $k+1$ vertices of degree $k-1$.  Now we can state Tverberg's lemma.

\begin{lem} 
\label{tver}
Let $G$ be a connected graph and let $k\DefinedAs\Delta(G)$.  If $k\ge 3$ and
$G$ is neither a $k$-tree nor $K_{k+1}$, then $G$ has a vertex $v$ of degree
$k$ such that no component of $G-v$ is a $k$-tree.
\label{ktrees-lemma}
\end{lem}

\begin{proof}
Suppose the lemma is false and let $G$ be a counterexample.  Now $G$ contains
a $T_k$, for otherwise any vertex $v$ of degree $k$ suffices.  Let $H$ be a copy
of $T_k$ with the minimum number of vertices of degree $k$ in $G$.  First,
suppose $H$ has only one vertex $v$ of degree $k$, and note that $v$ is a
cutvertex.  Since
$G$ is not a $k$-tree and $H$ is a $k$-tree, $G-H$ is not a $k$-tree. Since
neither $G-H$ nor $H-v$ is a $k$-tree, $v$ is the desired vertex, which is a
contradiction.

So instead $H$ has at least two vertices of degree $k$. Pick neighbors $v,w \in
V(H)$, where $v$ has degree $k$ and $w$ has degree as small as possible. When
$G-v$ is connected, let $A \DefinedAs G-v$; otherwise $G-v$ has two components,
$A$ and $B$, where $H-v \subseteq A$.   Since $G$ is a counterexample, either
$A$ or $B$ is a $k$-tree.

\begin{center}
\begin{figure}[h!tb]

\begin{tikzpicture}[line width=.03cm,scale=.65]
\tikzstyle{vert}=[draw,shape=circle,fill=black,minimum size=4pt, inner sep=0pt] 

\draw (3,0) arc (0:-170:3 and 2); 
\draw (3,0) arc (0:170:3 and 2); 
\draw (-2.45,0) arc (0:62:1 and .4); 
\draw (-2.45,0) arc (0:-62:1 and .4); 

\draw (-1.5,0) circle (1.3); 
\draw (-6.5,0) circle (1.3); 
\draw (-6.7,.2) node[] (B) {$B$}; 
\draw (-1.3,.2) node[] (H) {$H$}; 
\draw (1.3,.2) node[] (A) {$A$}; 

\draw (-2.8,0) node[vert,label={[label distance=-.10cm]150:{$v$}}] (v) {}; 
\draw (-1.0,-1.2) node[vert,label={[label distance=-.15cm]330:{$w$}}] (w) {};
\draw (-5.2,0) node[vert] (u) {}; 
\draw (u) -- (v);
\draw (v) -- (w);

\end{tikzpicture}

\caption{Components $A$ and $B$ of $G-v$ in the final case of
Lemma~\ref{ktrees-lemma}.} 
\end{figure}
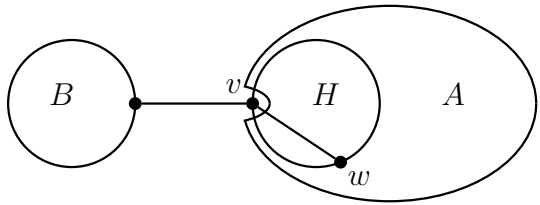
\end{center}

%
Suppose $B$ is a $k$-tree. Now we find a copy of $T_k$ with at most one vertex
of degree $k$, which contradicts the minimality of $H$.  
If $B=T_k$, then $B$ will do. 
Otherwise, we choose a leaf of $B$ with no vertex adjacent to $v$ in $G$.  
Thus $B$ is not a $k$-tree; so instead $A$ is a $k$-tree.  

Since $A$ is a $k$-tree, $w$ has degree at least $k-1$ in $A$ and hence at
least $k$ in $G$.  By our choice of $w$, every vertex in $H$ has degree $k$ in
$G$.  Hence, by the minimality of $H$, every vertex in a $T_k$ in $G$ has
degree $k$.  But every vertex of $A$ is in a $T_k$ and hence has degree $k$ in
$A$ unless it is a neighbor of $v$.  So, $A$ is a $k$-tree with at most $k$
vertices of degree $k-1$.  The only such $k$-tree is $T_k$, so $A = T_k$ and
$G=K_{k+1}$, a contradiction.
\end{proof}

\begin{proof}[Proof~\ref{ktrees} of Brooks' Theorem]
Suppose the theorem is false and choose a minimum counterexample $G$.  If $G$
is a $\Delta$-tree, then let $v$ be a cutvertex.
For each component $H$ of $G-v$, we can $\Delta$-color $H+v$, and then
permute the colors so the colorings agree on $v$.  So $G$ is not a
$\Delta$-tree.  Now apply Lemma
\ref{tver} with $k \DefinedAs \Delta$ recursively on components until all
components have maximum degree less than $k$; let $v_1, \ldots, v_r$ be the
vertices used by the lemma and let $I \DefinedAs \set{v_1, \ldots, v_r}$.
Note that $I$ is independent.  If $k = 3$, then $G-I$ consists of even cyles
and paths, so we can $2$-color it.  If $k \ge 4$, then by minimality of
$\card{G}$, we can $(k-1)$-color $G-I$.  Now we finish by coloring $I$ with a
new color, a contradiction.
\end{proof}

\section{Partitioned coloring}
\label{partitioned}
The \emph{coloring number} $\operatorname{col}(G)$ is defined by
$\operatorname{col}(G) \DefinedAs 1+\max_{H\subseteq G}\delta(H)$, where $H$
ranges over all subgraphs of $G$. 
Suppose we color a
graph $G$ as follows: delete a vertex $v$ of minimum degree, recursively
color $G-v$, and greedily color $v$.  This method shows that $\chi(G)\le
\operatorname{col}(G)$.
In 1976, Borodin \cite{borodin1976decomposition} proved the following generalization of Brooks' Theorem.

\begin{thm}
\label{BorodinDecomp}
Let $G$ be a graph not containing a $K_{\Delta + 1}$. If $s \ge 2$ and $r_1, \dots, r_s$ are
positive integers with $r_1 + \dots + r_s \ge \Delta \ge 3$, then $V(G)$ can be 
partitioned into sets $V_1, \dots, V_s$ such that $\Delta(G[V_i]) \leq r_i$ and
$\operatorname{col}(G[V_i]) \le r_i$ for all $i \in \irange{s}$.  
\label{borodinPartition}
\end{thm}

Brooks' Theorem follows from Theorem~\ref{borodinPartition} by taking $s =
\Delta$ and $r_1 = \dots = r_s = 1$.\footnote{
Section
\ref{DegreeChoosableGraphSection} gives another extension of
Brooks' Theorem (also due to Borodin \cite{borodin1977criterion}) that
classifies the graphs that are colorable when each vertex $v$ is allowed a
list of colors $L(v)$ and $|L(v)| = d(v)$.  
Borodin, Kostochka, and Toft~\cite{borodin2000variable} proved an
intruiging common generalization of these two, seemingly unrelated, results.
}
One way to prove Theorem \ref{BorodinDecomp} is to extend the idea of the
fundamental result of Lov\'asz \cite{lovasz1966decomposition} along the lines
of Catlin \cite{CatlinAnotherBound} and Bollob\'as and Manvel
\cite{bollobasManvel}, where a partition is repeatedly modified by moving
vertices from one part to another.  In his dissertation the second author gave
the following proof of Brooks' Theorem, where this dynamic process is
specialized and made static.  (The heart of this proof is
Lemma~\ref{ObstructionFree}, which is similar to the special case of
Theorem~\ref{borodinPartition} when $s=2$, $r_1=1$, and $r_2=\Delta-1$.)

Let $G$ be a graph.  A partition $P \DefinedAs (V_1, V_2)$ of $V(G)$
is \emph{normal} if it minimizes the value of $(\DeltaG - 1)
\size{V_1}+\size{V_2}$, where $\size{V_i}$ denotes $\card{E(G[V_i])}$.
%
Note that if $P$ is a normal partition, then
$\Delta(G[V_1]) \leq 1$ and $\Delta(G[V_2]) \leq \DeltaG - 1$, since if some
vertex $v$ has degree that is too high in its part, then we can move it to the
other part.  
The \emph{$P$-components} of $G$ are the components of $G[V_1]$ and $G[V_2]$. 
A $P$-component is an \emph{obstruction} if it is a $K_2$ in $G[V_1]$ or a
$K_{\Delta}$ in $G[V_2]$ or an odd cycle in $G[V_2]$ when $\Delta = 3$.  (Note
that if a $P$-component contains an obstruction, then the obstruction is the
whole $P$-component.)

A path $v_1\dots v_k$ is \emph{$P$-acceptable} (see
Figure~\ref{partition-pics}) if $v_1$ is in an obstruction and for all $i, j
\in \{1,\dots,k\}$, $v_i$ and $v_j$ are in different $P$-components.  A
$P$-acceptable path is \emph{maximal} if it is not contained in a larger
$P$-acceptable path.  This means that a $P$-acceptable path $v_1\dots v_k$ is
maximal if and only if every neighbor of $v_k$ is in the same $P$-component as
some vertex in the path.  Given a partition $P$, to \emph{move} a vertex $u$ is
to move it to the other part of $P$.  Note that if $P$ is normal and $u$ is in
an obstruction, then the partition formed by moving $u$ is again normal since
$u$ had maximum degree in its original part.  For a subgraph $H$ of $G$ and
vertex $u \in V(G)$, let $N_H(u) \DefinedAs N(u) \cap V(H)$.

\begin{lem}
\label{ObstructionFree}
Let $G$ be a graph with $\Delta \geq 3$.  If $G$ does not contain
$K_{\Delta + 1}$, then $V(G)$ has an obstruction-free normal partition.
\end{lem}
\begin{proof}
Suppose the lemma is false and let $G$ be a counterexample.  Among the normal
partitions of $G$ with the minimum number of obstructions, choose $P \DefinedAs
(V_1, V_2)$ and a maximal $P$-acceptable path $v_1\dots v_k$ so as to minimize
$k$.  Throughout the proof, we often move some vertex $u$ in an obstruction
$A$.  Since this destroys $A$, the minimality of $P$ implies that the move
creates a new obstruction, which must contain $u$.  So if $u$ has a neighbor
$w$, initially in the other part, then $w$ is in this new obstruction.
Finally, the new partition has the minimum number of obstructions.

Let $A$ and $B$ be the $P$-components containing $v_1$ and $v_k$ respectively.  
Let $X \DefinedAs N_A(v_k)$.  If $\card{X} = 0$, then moving $v_1$ creates a
new normal partition $P'$. 
Since $v_1$ is adjacent to $v_2$, the new obstruction contains $v_2$.
So $v_2v_3\dots v_k$ is a maximal $P'$-acceptable
path, violating the minimality of $k$.  Hence $\card{X} \geq 1$.

Pick any vertex $x \in X$, and form $P'$ by moving $x$.  The new obstruction
contains $v_k$ and hence all of $B$, since each obstruction is a whole
$P'$-component.  So $\set{x} \cup V(B)$ induces an obstruction.  Let $Y
\DefinedAs N_B(x)$; see Figure~\ref{partition-pics}.  Since $x\in X$ was arbitrary, the argument works for all
$z\in X$.  Since obstructions are regular, $N_B(z) = Y$ for all $z \in X$,
which implies that $X$ is joined to $Y$ in $G$.
Also since obstructions are regular, $\card{Y} = \delta(B) + 1$.  

\begin{center}
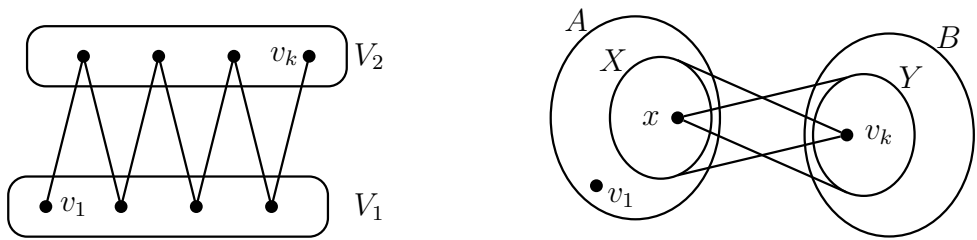
\begin{figure}[htb]

\begin{tikzpicture}[line width=.03cm, xscale=.5, yscale=1]
\tikzset{->-/.style={decoration={
  markings,
  mark=at position #1 with {\arrow{>}}},postaction={decorate}}}

\tikzstyle{vert}=[draw,shape=circle,fill=black,minimum size=4pt, inner sep=0pt] 

 \draw (1,0) node[vert,label={[label distance=-.05cm]0:$v_1$}] {} to (2,2);
 \draw (2,2) node[vert] {} to (3,0);
 \draw (3,0) node[vert] {} to (4,2);
 \draw (4,2) node[vert] {} to (5,0);
 \draw (5,0) node[vert] {} to (6,2);
 \draw (6,2) node[vert] {} to (7,0);
 \draw (7,0) node[vert] {} to 
       (8,2) node[vert,label={[label distance=-.1cm]180:$v_k$}] {};

\draw[rounded corners=.25cm] 
     (.5,1.6) -- (9,1.6) -- (9,2.4) -- (.5,2.4) -- cycle; 
\draw[rounded corners=.25cm] 
     (0,-.4) -- (8.5,-.4) -- (8.5,.4) -- (0,.4) -- cycle;

\draw (9.6,0) node {{$V_1$}};
\draw (9.6,2) node {{$V_2$}};
\end{tikzpicture}
\hspace{.7in}
\begin{tikzpicture}[line width=.03cm, scale = .45]
\tikzstyle{vert}=[draw,shape=circle,fill=black,minimum size=4pt, inner sep=0pt] 

\begin{scope}[yshift=.2in]
\draw \boundellipse{0,0}{1.5}{1.8}; 
\draw \boundellipse{-.75,0}{2.5}{3}; 
\draw (-1.4,1.7) node {$X$};
\draw (-2.5,2.9) node {$A$};
\draw (0.5,0) node[vert,label=left:{$x$}]    (x)  {};
\end{scope}

\draw \boundellipse{6,0}{1.5}{1.8}; 
\draw \boundellipse{6.75,0}{2.5}{3}; 
\draw (7.4,1.7) node {$Y$};
\draw (8.5,2.9) node {$B$};
\draw (5.5,0) node[vert,label=right:{$v_k$}] (vk) {};

\draw (x) to (5.5,1.7);
\draw (x) to (5.5,-1.7);

\draw (vk) to (.5,2.21);
\draw (vk) to (.5,-1.19);

\draw (-1.9,-1.5) node[vert,label={[label distance=-.2cm]273:$v_1$}] {}; 
\end{tikzpicture}

\caption{The left figure shows an example of a $P$-acceptable path.
The right figure shows sets $X$ and $Y$ constructed from
components $A$ and $B$.
}
\label{partition-pics}
\end{figure}
\end{center}

First, suppose $\card{X} \geq 2$.  Similar to above, form $P'$ from $P$ by
moving $x$ and $v_k$.  Now $\set{v_k} \cup V(A - x)$ induces an obstruction
($v_k$ is in a $P'$-component with $v_1$, since $|X-x|\ge 1$). 
Because obstructions are regular, $\card{N_{A-x}(v_k)} = \Delta(A)$
and hence $\card{X} \geq \Delta(A) + 1$.  Since $X$ and $Y$ are disjoint,
$\card{X \cup Y} \ge (\Delta(A) + 1) + (\delta(B) + 1) = \Delta(G) + 1$;  
the equality holds because $A$ and $B+x$ are obstructions.  
(If $\Delta > 3$, then $\Delta(A)+1=\card{A}$ and $\delta(B)+1=\card{B}$ and the
sizes of obstructions in distinct parts sums to $\Delta+2$; the case $\Delta=3$
is a little different, since now one obstruction is an odd cycle.)

Suppose $X$ is not
a clique and pick nonadjacent $x_1, x_2 \in X$. It is easy to check that moving
$x_1$, $v_k$, $x_2$, violates the normality of $P$.  Hence $X$ is a clique. 
Similarly, suppose $Y$ is not a clique and pick nonadjacent $y_1, y_2
\in Y$ and any $x' \in X - \set{x}$. Now moving $x$,~$y_1$,~$x'$,~$y_2$ 
again violates the normality of $P$.  Hence $Y$ is a clique.  But $X$ is
joined to $Y$, so $X \cup Y$ induces $K_{\Delta + 1}$ in $G$, a contradiction.

So instead $\card{X} = 1$.  Suppose $X \neq \set{v_1}$, and first suppose
$A$ is $K_2$.  Now moving $x$ creates another normal
partition $P'$ with the minimum number of obstructions.  In $P'$,
the path $v_kv_{k-1}\dots v_1$ is a maximal $P'$-acceptable path, since the
$P'$-components containing $v_2$ and $v_k$ contain all neighbors of $v_1$ in
that part.  Repeating the above argument using $P'$ in place of $P$ gets us to
the same point with $A$ not $K_2$.  Hence we may assume $A$ is not $K_2$.

Move each of $v_1, \dots, v_k$ in turn.  The obstruction $A$ is destroyed by
moving $v_1$, and for $1 \leq i < k$, the obstruction created by
moving $v_i$ is destroyed by moving $v_{i+1}$.  So after the moves, $v_k$
is in an obstruction.  The minimality of $k$ implies that
$\set{v_k} \cup V(A - v_1)$ induces an obstruction and hence $\card{X} \geq
2$, since $A$ is not $K_2$.  This contradicts $\card{X}=1$. 

Therefore $X = \set{v_1}$.  Now moving $v_1$ creates an obstruction
containing both $v_2$ and $v_k$,  so $k = 2$.
Since $v_1v_2$ is maximal, $v_2$ has no neighbor in the other part besides
$v_1$.  But now moving $v_1$ and $v_2$ creates a partition violating
the normality of $P$.
\end{proof}

\begin{proof}[Proof \ref{partitioned} of Brooks' Theorem.]
Suppose Brooks' Theorem is false and choose a counterexample $G$ minimizing
$\Delta$.  Clearly $\Delta \geq 3$.  
By Lemma \ref{ObstructionFree}, $V(G)$ has an
obstruction-free normal partition $(V_1, V_2)$.  
Note that $V_1$ is an independent set, since $\Delta(G[V_1])\le 1$ 
and $G[V_1]$ contains no $K_2$. 
Since $G[V_2]$ is obstruction-free,
the minimality of $\Delta$ gives $\chi(G[V_2]) \leq \Delta(G[V_2]) \le
\Delta-1$.  Using one more color on $V_1$ gives 
$\chi(G) \leq 1+\chi(G[V_2]) \le \Delta$, a contradiction.
\end{proof}

\section{Spanning trees with independent leaf sets}
\label{leafset}
In this section, we take a scenic route.  We combine two lemmas of
independent interest, to give an unexpected proof of Brooks' Theorem.
The union of two forests $F_1$ and $F_2$ (on the same vertex set) is
4-colorable, since we can color each vertex $v$ with a pair $(a_1,a_2)$, where
$a_i$ is the color of $v$ in a proper 2-coloring of $F_i$.
A \emph{star forest} is a disjoint union of stars.
Sauer conjectured that the union of a forest and a star forest is always
$3$-colorable, and Stiebitz \cite{forestPlusStars} verified this conjecture.
The main component of his proof is a lemma that allows for extending
a $k$-coloring of an induced subgraph to the whole graph when it has a
spanning forest with certain properties.  

B\"{o}hme et al.~\cite{independencyTree} classified the graphs with a
spanning tree whose leaves form an independent set; by combining this result
with Stiebitz's coloring lemma, they gave an alternative proof of Brooks'
Theorem.
In this section, we prove both the lemma of Stiebitz and that of B\"{o}hme et
al., as well as show how they easily yield Brooks' Theorem.

We need two definitions.
An \emph{independency tree} is a spanning tree in which the leaves form an
independent set.  Let $v_1,\ldots,v_n$ be an order of the vertices of a graph
$G$; call it
$\sigma$.  We define a \emph{depth-first-search tree}, or \emph{DFS tree}, with
respect to $\sigma$ as follows.  We iteratively grow a tree. 
At each step $i$, we have a tree $T_i$ and an \emph{active vertex} $x_i$.
Let $T_1\DefinedAs v_1$ and $x_1\DefinedAs v_1$.  We grow the tree as follows.
If $x_i$ has a neighbor not in $T_i$, then choose the first such neighbor
$w$, with respect to $\sigma$.  
Now let $T_{i+1}\DefinedAs T_i+x_iw$ and $x_{i+1} \DefinedAs w$.  
If $x_i$ has no neighbor outside the tree, then let $w$ be the neighbor of
$x_i$ on a path in $T_i$ to $v_1$.  Now let $T_{i+1}\DefinedAs T_i$ and
$x_{i+1}\DefinedAs w$.  
This algorithm terminates only when $x_k=v_1$ and all neighbors of $v_1$ are
in the tree.  It is easy to check that when this happens $T_k$ is a
spanning tree.  A \emph{DFS independency tree} is both a DFS tree and an
independency tree.
We begin with the following elegant lemma, from B\"{o}hme et
al.~\cite{independencyTree}; about 25 years earlier, Dirac and
Thomassen~\cite{DiracT73} proved a variation containing (1) and (4), as well as
other equivalent conditions (but not (2) or (3)).

\begin{lem}
For a connected graph $G$, the following four conditions are equivalent.
\begin{enumerate}
\item $G$ is $C_n$, $K_n$, or $K_{n/2,n/2}$ for $n$ even.
\item $G$ has no independency tree.
\item $G$ has no DFS independency tree.
\item $G$ has a Hamiltonian path, and every Hamiltonian path is contained in
a Hamiltonian cycle.
\end{enumerate}
\label{independencytree}
\end{lem}
\begin{proof}
(1) $\implies$ (2): If $G$ is $C_n$ or $K_n$, then any spanning tree has 
all leaves pairwise adjacent, so $G$ has no independency tree.  So Let $G$
be $K_{n/2,n/2}$.  If all leaves are in the same part, then each vertex in the
other part has degree at least 2. This requires at least $n$ edges, which is
too many for an $n$-vertex tree; so we get a contradiction.
\smallskip

(2) $\implies$ (3): This implication is immediate from the definitions.

(3) $\implies$ (4): We prove the contrapositive.  
If $G$ has a Hamiltonian path $P$ that is not contained in a Hamiltonian cycle,
then $P$ is an independency tree.  Considering the vertices in the order in
which they appear in $P$ shows that $P$ is a DFS independency tree.  So
suppose instead that $G$ has no Hamiltonian path.  Let $P=v_1\ldots v_k$ be a
maximum path in $G$.  The maximality of $P$ implies that $v_1$ and $v_k$ are
nonadjacent to each $u\in V\setminus V(P)$. Similarly $v_1$ and $v_k$ are
nonadjacent (if not, let $u$ be a neighbor of some $v_i$ not on $P$; now the
path $v_{i+1}\ldots v_kv_1\ldots v_iu$ is longer than $P$).  Now consider a DFS
tree $T$ that begins with $v_1,\ldots,v_k$.  Since it has an independent leaf
set, $T$ is a DFS independency tree.

\input{cycle-figs.tex}
(4) $\implies$ (1):
Suppose that $G$ has a Hamiltonian path $P$,  
and $P$ is contained in a Hamiltonian cycle $C$, where $C=v_1\ldots v_n$.  If
$C$ has no chords, then $G=C_n$, and (1) holds.  So assume $C$ has a chord. 
The \emph{length} of a chord $v_iv_j$ is $\min(|i-j|,n-|i-j|)$. 

The left part of Figure~6 shows that if $G$ contains a chord of $C$ of a given
length, then G contains all chords of $C$ of that length; the central part
shows that if $G$ contains a chord of $C$ of odd (even) length, then it
contains all chords of $C$ of odd (even) length.  

Suppose $n$ is odd. If $G$ has a chord, then it has either $v_nv_{(n+3)/2}$ or
$v_nv_{(n-1)/2}$ (since one chord has odd length and the other even length).
By the middle of Figure~6, the presence of one of these chords implies the
presence of the other.  So $G$ has all chords and $G=K_n$.
Suppose instead that $n$ is even. If $G$ has an even chord, then $G$ has a
chord of length 2. 
For any chord $v_nv_k$ of length at least 3, the right part of Figure~6 shows
that $G$ also contains $v_nv_{k+1}$; hence $G$ contains chords of both
parities, so $G=K_n$.
In the final case, with $n$ even and no even chord, all odd chords exist, so
$G=K_{n/2,n/2}$.
\end{proof}

%
%

%

Next we state Stiebitz's lemma for extending a $k$-coloring of a subgraph to the
whole graph.

\begin{lem}
\label{stiebitz}
Let $H$ be an induced subgraph of a graph $G$ with $\chi(H) \le k$ for some $k \ge 3$.  Then $\chi(G) \le k$ if $G$ has a spanning forest $F$ where
\begin{enumerate}
\item for each component $C$ of $H$, $F\brackets{V(C)}$ is a tree; and
\label{cond1}
\item $d_G(v) \le d_F(v) + k-2$ for every $v \in V(G-H)$.
\label{cond2}
\end{enumerate}
\end{lem}

Before proving Lemma \ref{stiebitz}, we 
use Lemmas~\ref{independencytree} and \ref{stiebitz}
to give a short proof of Brooks' Theorem.
\begin{proof}[Proof \ref{leafset} of Brooks' Theorem]
It suffices to consider connected graphs.
If a graph $G$ is $C_n$, $K_n$, or $K_{n/2,n/2}$ for $n$ even, then $\chi(G)$
satisfies the desired bound.  So suppose $G$ is none of these graphs.  By
Lemma~\ref{independencytree}, $G$ has an independency tree $T$.  Let
$I$ denote the independent leaf set of $T$.  To apply Lemma~\ref{stiebitz},
let $H\DefinedAs G[I]$, let $F\DefinedAs T$, and let
$k\DefinedAs \Delta$.  Clearly $\chi(H)=1\le k$.  Each component $C$ of $H$ is
an isolated vertex, which is a tree, so (1) holds.  Finally, each vertex $v\in
V(G-H)$ is a nonleaf in $F$, so $d_F(v)\ge 2$.  Thus $d_G(v)\le k \le d_F(v) +
(k-2)$.  So Lemma~\ref{stiebitz} implies that $\chi(G)\le k = \Delta$.
\end{proof}

The proof of Lemma~\ref{stiebitz} yields an algorithm which extends the
coloring of $H$ to one of $G$, by adding vertices to $H$ one at a time (this is
a rough approximation; we give more precise details below).
The proof is by contradiction and it relies on fives claims; any vertex $v$ 
violating a claim allows us to make progress in extending the coloring.
One obvious route is to color $v$ immediately and add it to $H$.  When this is
not possible, we form a smaller graph $G'$ from $G-v$ by identifying some
vertices, and we color $G'$ recursively; afterwards, we color $v$ greedily. 
Clearly $G'$ (and $F'$ and $H'$) must satisfy the hypotheses of the lemma.  We
also must ensure that the resulting coloring of $G-v$ uses at most $k-1$ colors
on neighbors of $v$, so that we can extend the coloring to $v$.

\begin{proof}[Proof of Lemma~\ref{stiebitz}.]
For any graphs $U$ and $W$, we write $U-W$ for the subgraph of $U$ induced by
$V(U)\setminus V(W)$.
If $uv\in E(F)$, then $u$ is an \emph{$F$-neighbor} of $v$, and $u$ and $v$ are
\emph{$F$-adjacent}.
Suppose the lemma is false and choose a counterexample pair $G, H$ minimizing
$|G-H|$.
%
Note that each vertex $v$ in $G-H$ must have a neighbor in $H$, since 
otherwise we can add $v$ to $H$.  Thus $\card{H}\ge 1$.

\textbf{Claim 1.} \textit{If there exists $v \in V(G-H)$ adjacent to
components $A_1, \ldots, A_s$ of $H$ with $d_G(v) \le s + k - 2$, then there
exist $i$ and $j$, with $i \ne j$, and a path in $F-v$ from $A_i$ to $A_j$.}

Suppose not and choose such a $v \in V(G-H)$.  
We will find a $k$-coloring of $G$.
For each $i \in \irange{s}$, let $z_i$ be a neighbor of $v$ in $A_i$.  Form
$G'$, $F'$, $H'$ from $G$, $F$, $H$ (repectively) by deleting $v$ and
identifying all $z_i$ as a single new vertex $z$.  
Now $\chi(H') \le k$, since by permuting colors in each component we can get
a $k$-coloring of $H$ where all the $z_i$ use the same color.  Also, $F'$ is a
spanning forest in $G'$ since we are assuming there is no path in $F-v$ from
$A_i$ to $A_j$ whenever $i \ne j$.  It is easy to check that Conditions 
\hyperref[cond1]{(1)} and \hyperref[cond2]{(2)} hold for $G', F', H'$.  Now
$|G'-H'|<|G-H|$, so by minimality of $|G-H|$,
we have a $k$-coloring of $G'$.  This gives a $k$-coloring of $G-v$ where $z_1,
\ldots, z_s$ all get the same color.  So $v$ has at most $d_G(v) - (s-1) \le
k-1$ colors used on its neighborhood, leaving a color free to finish the
$k$-coloring on $G$, a contradiction.

\begin{center}
\begin{figure}[h!tb]

\begin{tikzpicture}[line width=.03cm,scale=.5] 
\tikzstyle{vert}=[draw,shape=circle,fill=black,minimum size=4pt, inner sep=0pt] 
\tikzstyle{small}=[draw,shape=circle,fill=black,minimum size=2pt, inner sep=0pt] 

\draw \boundellipse{-4,0}{3.5}{2.5}; 
\draw (-.9,0) node[vert,label={[label distance=-.05cm]180:{$v$}}] (v) {}; 
\draw (-4.5,.5) node (GminusH) {$G-H$}; 

\draw \boundellipse{2.5,1.65}{1.5}{.9}; 
\draw (4.7,1.7) node (A1) {$A_1$}; 
\draw (1.5,1.65) node[vert,label={[label distance=-.1cm]0:{$z_1$}}] (z1) {}; 

\draw \boundellipse{2.5,-1.65}{1.5}{.9}; 
\draw (4.7,-1.7) node (As) {$A_s$}; 
\draw (1.5,-1.65) node[vert,label={[label distance=-.1cm]0:{$z_s$}}] (zs) {}; 

\draw (z1) -- (v) -- (zs);
\draw (2.0,.4) node[small] () {};
\draw (2.0,0) node[small] () {};
\draw (2.0,-.4) node[small] () {};
\end{tikzpicture}
%
\hspace{.25in}
\begin{tikzpicture}[line width=.03cm,scale=.35] 
\tikzstyle{vert}=[draw,shape=circle,fill=black,minimum size=4pt, inner sep=0pt] 

\draw \boundellipse{-6,.5}{1.3}{1};
\draw (-6.3,.2) node[vert] (u1) {}; 
\draw \boundellipse{-2,.5}{1.3}{1};
\draw (-2,.2) node[vert] (u2) {}; 

\draw[rounded corners=.7cm] (-8,-2) -- (-8,-6) -- (0,-6) -- (0,-2) -- cycle;
\draw (-6,-3) node[vert,label={[label distance=-.05cm]360:{$v$}}] (v) {}; 
\draw (-5.5,-5) node[vert] (v2) {}; 
\draw (-2.5,-5) node[vert] (v3) {}; 
\draw (-2,-3) node[vert] (v4) {}; 
\draw (u1) -- (v) -- (v2) -- (v3) -- (v4) (u2) -- (v);

\begin{scope}[xshift=10cm] 

\draw \boundellipse{-6,.5}{1.3}{1};
\draw (-6.7,.2) node[vert] (u3) {}; 
\draw (-6.0,.4) node[vert] (u4) {}; 
\draw (-5.3,.2) node[vert] (u5) {}; 
\draw \boundellipse{-2,.5}{1.3}{1};
\draw (-2,.2) node[vert] (u6) {}; 

\draw[rounded corners=.7cm] (-8,-2) -- (-8,-6) -- (0,-6) -- (0,-2) -- cycle;

\draw (-5.0,-3) node[vert] (v5) {}; 
\draw (-2.5,-3) node[vert] (v6) {}; 
\draw (v4) -- (u3) -- (u4) -- (u5) -- (v5) -- (v6) -- (u6);
\end{scope}

\draw[rounded corners=.2cm] (-8.5,-1.5) -- (-8.5,-6.5) -- (10.5,-6.5) --
(10.5,-1.5) -- cycle; 
\draw[rounded corners=.2cm] (-8.5,-1) -- (-8.5,2) -- (10.5,2) --
(10.5,-1) -- cycle; 

\draw (12.5,-4) node[] (GminusH) {$G-H$}; 
\draw (11.5,.5) node[] (H) {$H$}; 
\draw (-7,-4.5) node[] (A) {$A$}; 
\end{tikzpicture}

\caption{The left figure shows Claim 1. The right figure 
shows Claim 3.
}
\label{lemma7-figs}
\end{figure}
\end{center}

\textbf{Claim 2.} \textit{Every leaf of $F$ is in $H$ and every vertex not in
$H$ has an $F$-neighbor not in $H$.}
We can rewrite this formally: $d_F(v) \ge 2$ and $d_{F-H}(v) \ge 1$ for all $v
\in V(G-H)$.
Applying Claim 1 with $s=1$ implies $d_G(v) \ge k$. 
Now 

{\CondTwo} gives $d_F(v) \ge d_G(v) + 2 - k \ge 2$. Suppose
$d_{F-H}(v) = 0$ for some $v \in V(G-H)$.   Since $F$ is a forest, {\CondOne} 
implies that all $F$-neighbors of $v$ must be in different components of $H$.
Moreover there can be no path between two of these components in $F-v$.  
{\CondTwo} gives $d_G(v)\le d_F(v)+k-2$, so
applying Claim 1 with $s = d_F(v)$ gives a contradiction.
Thus $d_{F-H}(v) \ge 1$ for all $v\in V(G-H)$.

\textbf{Claim 3.} 
\textit{There exists $v$ in $G-H$ with $d_{F-H}(v) = 1$ such that every
component of $H$ that is $F$-adjacent to $v$ is not $F$-adjacent to any other
vertex in $G-H$.}

Form a bipartite graph $F'$ from $F$ by contracting each component of $H$ and
each component of $F-H$ to a single vertex.  Since $F$ is a forest, {\CondOne}
implies that $F'$ is also a forest.  So some vertex contracted from a component
$A$ of $F-H$ has at most one neighbor of degree at least 2; say this neighbor
is contracted from $B$, where $B\subseteq (F\cap H)$.  (If not, then we can
walk between components of $H$ and $F-H$ until we get a cycle in $F$.)
%
Let $v$ be a leaf of $A$ that is not $F$-adjacent to $B$; this gives
$d_{F-H}(v) = d_{A}(v)\le 1$.  Claim 2 gives $d_{F-H}(v)\ge 1$, so in fact
$d_{F-H}(v)=1$ as desired.

\textbf{Claim 4.} \textit{If the $v$ in Claim 3 is adjacent to a component of
$H$, then it is $F$-adjacent to that component.}

Let $A_1, \ldots, A_r$ be the components of $H$ that are $F$-adjacent to $v$,
where $r = d_F(v) - 1$.  Suppose there is another component $A_{r+1}$ of $H$
that is adjacent to $v$.  Since no vertex of $G-H$ besides $v$ is $F$-adjacent
to any of $A_1, \ldots, A_r$, there can be no $F$-path in $F-v$ between any
pair among $A_1, \ldots, A_r, A_{r+1}$.  Now the contrapositive of Claim 1
implies that $d_G(v) > (r + 1) + k - 2 = d_F(v) + k - 2$; this inequality
contradicts {\CondTwo}.

\textbf{Claim 5.} \textit{The lemma holds.}

Let $H' \DefinedAs G[V(H) \cup \set{v}]$, with $v$ as in Claims 3 and 4.  
By Claim 4, {\CondOne} of the
hypotheses holds for $H'$. {\CondTwo} clearly holds and $F$ is still a forest. Also, by permuting colors in the components we can get a $k$-coloring of $H$ where all $F$-neighbors of $v$ get the same color.  Hence $v$ has at most $d_H(v) - (d_F(v) - 2) \le d_G(v) - 1 - (d_F(v) - 2) = d_G(v) - d_F(v) + 1 \le k-1$ colors on its neighborhood.  Hence $H'$ is $k$-colorable. But then, by minimality of $|G-H|$, $G$ is $k$-colorable, a contradiction.
\end{proof}

\section{Kernel perfection}
\label{kernel}

In the late 1970s, Vizing \cite{Vizing76} and, independently, 
\erdos,
Rubin, and Taylor~\cite{ErdosRT79} introduced the notion of list coloring,
which is the subject of Sections~\ref{kernel}
and~\ref{DegreeChoosableGraphSection}.
An \emph{$f$-list assignment} gives to each vertex $v$ a list $L(v)$ of $f(v)$
allowable colors.  A \emph{proper $L$-coloring} is a proper coloring where each
vertex gets a color from its list.  A graph $G$ is \emph{$f$-choosable} (or
\emph{$f$-list colorable}) if it has a proper $L$-coloring for each $f$-list
assignment $L$.  We are particularly interested in two cases of $f$.  
If $G$ is $f$-choosable and $f$ is constant, say $f(v)=k$ for all $v$, then $G$
is $k$-choosable.  
The minimum $k$ such that $G$ is $k$-choosable is its \emph{choice number}
$\chil$.\footnote{\erdos, Rubin, and Taylor noted that bipartite graphs can have
arbitrarily large choice number.  Let $m={2k-1  \choose k}$ and let $G=K_{m,m}$.
If we assign to the vertices of each part the distinct $k$-subsets of
$\irange{2k-1}$, then $G$ has no good coloring.  Thus $\chil(G) > k$.}
If $f(v)=d(v)$ for all $v$ and $G$ is $f$-choosable, then $G$ is
\emph{degree-choosable}.  All our remaining proofs show that Brooks' Theorem is
true even for list coloring\footnote{Formally: If $G$ is a connected graph other
than an odd cycle or a clique, then $\chil(G)\le \Delta$.}, which is a stronger result.

In 2009,
Schauz~\cite{Schauz09} and Zhu~\cite{Zhu09online}  
introduced \emph{online list coloring}.
In this variation, list sizes are fixed (each vertex
$v$ gets $f(v)$ colors), but the lists themselves are provided
online by an adversary.  In round 1, the adversary reveals the set of 
vertices whose lists contain color 1.  The algorithm then uses
color 1 on some independent subset of these vertices (and cannot change this
set later).  In each subsequent round $k$, the adversary reveals 
the subset of uncolored vertices with lists containing $k$. 
Again the algorithm chooses an independent subset of these vertices 
to use color $k$.
The algorithm wins if it succeeds in coloring all vertices.  And the
adversary wins if it reveals a vertex $v$ on each of $f(v)$ rounds, but the
algorithm never colors it.  A graph is \emph{online $k$-list colorable} (or
\emph{$k$-paintable}) if some algorithm can win whenever $f(v)=k$ for all $v$. 
The minimum $k$ such that a graph $G$ is online $k$-list colorable is its
\emph{online list chromatic number}, denoted $\chiol$, (or \emph{paint
number}).\footnote{At first glance, an adversary assigning lists to vertices
seems to have much more power in the context of online list coloring than in
list coloring.  However, in practice the choice number and paint number are
often equal.  In fact, it is unknown \cite{Westetal-paint14} whether there
exists a graph with $\chiol > \chil+1$.}

Any proof of the following lemma yields a proof of Brooks' Theorem for
coloring, list coloring, and even online list coloring using the kernel ideas
below (we write $\alpha(G)$ to denote the maximum size of an independent set,
or stable set, in $G$).  Note that this lemma follows immediately from Brooks'
Theorem for coloring and also from Brooks' Theorem for fractional coloring
(defined below).  Surprisingly, all known short proofs of Lemma~\ref{ind-ratio}
rely on some version of Brooks' Theorem.

\begin{lem}\label{BrooksIndependentSet}
If $G$ is a graph not containing $K_{\DeltaG + 1}$, then $\alpha(G) \ge
\frac{|G|}{\DeltaG}$.
\label{ind-ratio}
\end{lem}

It's natural to ask when we can improve the bound in Lemma~\ref{ind-ratio}.
When $G$ can be partitioned into disjoint copies of $K_{\Delta}$, clearly we
cannot.  But Albertson, Bollobas, and Tucker~\cite{AlbertsonBT76} did improve
the bound when $G$ is connected and $K_{\Delta}$-free, except when $G$ is one
of the two
exceptional graphs shown in Figure~\ref{exceptions-fig}.  
Along these lines, Fajtlowicz~\cite{Fajtlowicz78} proved that every graph
$G$ satisfies $\alpha(G) \ge \frac{2|G|}{\omega+\Delta+1}$.  In general, this
bound is weaker than Lemma~\ref{ind-ratio}, but when $\omega \le
\Delta-2$ it is stronger.

\begin{figure}[!ht]
\begin{center}
\begin{tikzpicture}[scale=.4,minimum size = 2pt, inner sep=0pt]
   \GraphInit[vstyle=Hasse]

\tikzset{VertexStyle/.style = {shape = circle,fill = black,minimum size = 5pt,inner sep=0pt}}

\begin{scope}[rotate=22.5]
\grCycle[prefix=a,RA=4.15]{8}
\EdgeInGraphMod{a}{8}{2}
\end{scope}

\begin{scope}[xshift = 5in, rotate=18]
\grCycle[prefix=a,RA=3.0]{5}
\grCycle[prefix=b,RA=5.0]{5}
\EdgeMod{a}{b}{5}{1}
\EdgeMod{a}{b}{5}{4}
\EdgeIdentity{a}{b}{5}
\end{scope}
\end{tikzpicture}

\end{center}
\caption{The only two connected $K_{\Delta}$-free graphs where
$\alpha = \frac{\card{G}}{\Delta}$.}
\label{exceptions-fig}
\end{figure}

A closely related problem is \emph{fractional coloring}, where independent
sets are assigned nonnegative weights so that for each vertex $v$ the sum of
the weights on the sets containing $v$ is 1.  In a \emph{fractional
$k$-coloring}, the sum of all weights on the independent sets is $k$; the
\emph{fractional chromatic number} is the minimum value $k$ allowing a
fractional $k$-coloring.  (A
standard vertex coloring is the special case when the weight on each set is 0
or 1.)  King, Lu, and Peng~\cite{KingLP12} strengthened the result of Albertson
et al.\ by showing that every connected $K_{\Delta}$-free graph with $\Delta\ge
4$ (except for the two graphs in Figure~\ref{exceptions-fig}) has fractional
chromatic number at most $\Delta-\frac{2}{67}$; this result was further
strengthened by Edwards and King~\cite{EdwardsK13,EdwardsK14+}.  
Recently Dvo\v{r}\'{a}k, Sereni, and Volec~\cite{DvorakSV14+} proved
fascinating results on fractional coloring of triangle-free cubic graphs. 
Specifically, they proved that all such graphs have fractional chromatic number
at most $\frac{14}5$, which is best possible.

Kostochka and Yancey~\cite{kostochkayancey2012ore} gave a simple, yet
powerful, application of the Kernel Lemma to a list coloring problem.  A
\emph{kernel} in a digraph $D$ is an independent set $I \subseteq V(D)$ such
that each vertex in $V(D) - I$ has an edge into $I$.  A digraph in which every
induced subdigraph has a kernel is \emph{kernel-perfect}, and
kernel-perfect orientations can be very useful for list coloring; 
 Alon and Tarsi~\cite[Remark 2.4]{AlonT92} attribute this result to Bondy,
Boppana, and Siegel (here $d^+(v)$ denotes the out-degree of $v$ in $D$).

\begin{KernelLemma}
Let $G$ be a graph and $\func{f}{V(G)}{\IN}$. If $G$ has a kernel-perfect orientation such that $f(v) \geq d^+(v) + 1$ for each $v \in V(G)$, then $G$ is $f$-choosable.
\label{kernel-lemma}
\end{KernelLemma}

The proof of the Kernel Lemma is by induction on the total number of colors in
the union of all the lists.  For an arbitrary color $c$, let $H$ be the
subdigraph induced by the vertices with $c$ in their lists.  We use $c$ on the
vertices of some kernel of $H$, and invoke the induction hypothesis to color
the remaining uncolored subgraph.  The same proof works for online list
coloring, since we simply choose $c$ to be color 1.

All bipartite digraphs are kernel-perfect, and the next lemma 
\cite{kostochkayancey2012ore} generalizes this fact.

\begin{lem}\label{KernelPerfect}
Let $A$ be an independent set in a graph $G$ and let $B \DefinedAs V(G) - A$. 
Any digraph $D$ created from $G$ by replacing each edge in $G[B]$ by a pair of
opposite arcs and orienting the edges between $A$ and $B$ arbitrarily is
kernel-perfect.
\end{lem}
\begin{proof}
Let $G$ be a minimum counterexample, and let $D$ be a digraph created from $G$
that is not kernel-perfect.  To get a contradiction it suffices to construct a
kernel in $D$, since each subdigraph has a kernel by minimality.  Either $A$ is
a kernel or there is some $v \in B$ which has no outneighbors in $A$.  In the
latter case, each neighbor of $v$ in $G$ has an inedge to $v$, so a kernel in
$D - v - N(v)$ together with $v$ is a kernel in $D$.
\end{proof}

Now we can show that Lemma \ref{BrooksIndependentSet} implies Brooks' Theorem. 
For simplicity we only prove this for list coloring, but a minor modification
of the proof works for online list coloring.
This proof is a special case of a result of the second author and
Kierstead~\cite{orevizing}.

\begin{center}
\begin{figure}[h!tb]

\begin{tikzpicture}[line width=.03cm]
\tikzset{->-/.style={decoration={
  markings,
  mark=at position #1 with {\arrow{>}}},postaction={decorate}}}

\tikzstyle{vert}=[draw,shape=circle,fill=black,minimum size=4pt, inner sep=0pt] 

 \draw[->-=.7] (1,0) node[vert] {} to (1,2);
 \draw[->-=.3] (1,2) node[vert] {} to (2,0);
 \draw[->-=.7] (2,0) node[vert] {} to (2,2);
 \draw[->-=.3] (2,2) node[vert] {} to (1,0);

 \draw[->-=.7] (3,0) node[vert] {} to (3,2);
 \draw[->-=.3] (3,2) node[vert] {} to (4,0);
 \draw[->-=.7] (4,0) node[vert] {} to (4,2);
 \draw[->-=.3] (4,2) node[vert] {} to (5,0);
 \draw[->-=.7] (5,0) node[vert] {} to (5,2);
 \draw[->-=.3] (5,2) node[vert] {} to (3,0);

 \draw[->-=.54] (1,2) to [bend left=90] (3,2);
 \draw[->-=.55] (3,2) to [bend right=45] (1,2);

 \draw[->-=.53] (2,2) to [bend left=55] (5,2);
 \draw[->-=.53] (5,2) to [bend right=35] (2,2);

\draw[rounded corners] (.8,1.8) -- (5.2,1.8) -- (5.2,2.2) -- (.8,2.2) -- cycle;
\draw[rounded corners] (.8,-.2) -- (5.2,-.2) -- (5.2,.2) -- (.8,.2) -- cycle;

\draw (5.6,0) node {\tiny{$A_H$}};
\draw (5.6,2) node {\tiny{$B_H$}};
\end{tikzpicture}

\caption{The bipartite graph $H$ with parts $A_H$ and $B_H$.
}
\label{kernel-pic}
\end{figure}
\end{center}
\begin{thm}
\label{KernelPerfectionBrooks}
Every graph satisfies $\chi_l \le \max\set{3, \omega, \Delta}$.
\end{thm}
\begin{proof}[Proof \ref{kernel} of Brooks' Theorem.]
Suppose the theorem is false and let $G$ be a minimum counterexample.  
The minimality of $G$ implies
$\chi_l(G - v) \leq \Delta$ for all $v \in V(G)$. 
So $G$ is $\Delta$-regular. 
Lemma \ref{BrooksIndependentSet} implies $\alpha(G) \ge \frac{|G|}{\Delta}$.
 Let $A$ be a maximum independent set in $G$ and let $B \DefinedAs V(G) - A$.
For each subgraph $H$, let $A_H\DefinedAs A\cap V(H)$ and $B_H\DefinedAs B\cap
V(H)$.
The number of edges between $A$ and $B$ is $\alpha(G)\Delta$, which is at least
$|G|$.  So 
there exists a nonempty induced subgraph $H$ of $G$ with at least $|H|$ edges
between $A_H$ and $B_H$, since $H\DefinedAs G$ is one example. 
Pick such an $H$ minimizing $|H|$.  
See Figure~\ref{kernel-pic}.

For all $v \in V(H)$, let $d_v$ be the number of edges incident to $v$ between
$A _H$ and $B_H$. 
We show that $d_v=2$ for all $v$. 
If $d_v\le 1$ for some $v$, then $H-v$ violates the minimality of $H$.  
The same is true if $d_v < d_w$ for some $v$ and $w$, as we now show.
Let $k$ be the minimim degree in $H$, and choose $v$ and $w$ with $d_v=k$ and
$d_w > d_v$. 
Let $\size{H}$ denote the number of edges in $H$.
Since $d_w > k$, we have
$2\size{H} \ge k|H| + 1$.
Deleting $v$ gives 
$2\size{H - v}  \ge k\card{H} + 1 - 2k = k(\card{H} - 1) + (1 - k) = 2(\card{H}
- 1) + (1 - k) + (k - 2)(\card{H} -  1)$.  
When $k \ge 3$ the final term is at least $k$, since $\card{H} \ge k + 1$; 
so $\size{H-v} \ge \card{H-v}$, a contradiction. 
If instead $k=2$, then $2\size{H-v} \ge 2(\card{H}-1)-1$.  Since the left side
is even, we conclude $2\size{H-v} \ge 2(\card{H}-1)$, again reaching a
contradiction.

So $d_v = d_w$ for all $v$ and $w$; we 
call this common value $d$, and note that $d\ge 2$.  
Now there are $(d/2)|H|$ edges between $A_H$ and $B_H$.  It is
easy to check that $(d/2)|H| - d \ge |H|
- 1$ for $d > 2$, so the minimality of $|H|$ shows that $d=2$.  Thus the
edges between $A_H$ and $B_H$ induce a disjoint union of cycles.

Create a digraph $D$ from $H$ by replacing each edge in $H[B_H]$ by a
pair of opposite arcs and orienting the edges between $A_H$ and $B_H$
consistently along the cycles.  
By Lemma \ref{KernelPerfect}, $D$ is kernel-perfect.
Since $d^+(v) \le d(v) - 1$ for each $v \in V(H)$,
the Kernel Lemma shows that $H$ is $f$-choosable where $f(v)
\DefinedAs d(v)$ for all $v \in V(H)$.  Now given any
$\Delta$-list-assignment on $G$, we can color $G-H$ from its lists by
minimality of $|G|$, and then extend the coloring to $H$, which is a
contradiction.
\end{proof}

\section{Degree-choosable graphs}
\label{DegreeChoosableGraphSection}
The option to give different vertices lists of different sizes allows us to
refine Theorem~\ref{KernelPerfectionBrooks}.
A graph $G$ is \emph{degree-choosable} if it is $f$-choosable where $f(v)
\DefinedAs d(v)$ for all $v \in V(G)$. The degree-choosable graphs were classified by Borodin \cite{borodin1977criterion} and 
independently by \erdos, Rubin, and Taylor \cite{ErdosRT79}. 
%
A graph is a \emph{Gallai tree} if each block is a complete graph or an odd
cycle.

\begin{center}
\begin{figure}[h!tb]
\begin{tikzpicture}[scale=.4,minimum size = 2pt, inner sep=0pt, line width=.8pt]
   \GraphInit[vstyle=Hasse]
\tikzstyle{vert}=[draw,shape=circle,fill=black,minimum size=4.3pt,inner sep=0pt] 
\tikzset{VertexStyle/.style = {shape = circle,fill = black,minimum size = 5pt,inner sep=0pt}}

\begin{scope}[rotate=38.57]
\grCycle[prefix=a,RA=4.15]{7}
\end{scope}

\begin{scope}[xshift=5.2cm,yshift=6.8cm,rotate=-17] 
\grCycle[prefix=b,RA=2.05]{3}

\begin{scope}[rotate=17] 
\begin{scope}[xshift=3.27cm,yshift=.95cm]
\begin{scope}[rotate=50]
\grCycle[prefix=c,RA=2.0]{4}
\EdgeInGraphMod{c}{4}{2}
\end{scope}
\end{scope}
\end{scope}
\end{scope}

\begin{scope}[xshift=15cm,yshift=8.1cm] 
\begin{scope}[rotate=10]
\grCycle[prefix=d,RA=2.2]{5}
\EdgeInGraphMod{d}{5}{2}
\end{scope}
\end{scope}

\begin{scope}[xshift=11.8cm, yshift=2cm] 
\begin{scope}[rotate=18] 
\grCycle[prefix=e,RA=2.1]{5}
\end{scope}
\end{scope}

\begin{scope}[xshift=21.2cm,yshift=5.8cm,rotate=-50] 
\grCycle[prefix=f,RA=2.0]{3}
\end{scope}

\begin{scope}[xshift=24.45cm,yshift=4.1cm,rotate=-65] 
\grCycle[prefix=g,RA=2.0]{3}
\end{scope}

\draw (a1) -- (b2) (c3) -- (d3) -- (e1) -- (c3)
(b1) -- ++(-1.75,2.5) node[vert] (b4) {} (b1) -- ++(.6,3.1) node[vert] (b5) {}
(d4) -- (f2);

\draw[rotate around={-65:(g0)}] 
(g0) -- ++ (2.5,1) node[vert] (h1) {}
(g0) -- ++ (2.5,0) node[vert] (h2) {}
(g0) -- ++ (2.5,-1) node[vert] (h3) {};
\end{tikzpicture}

\caption{A Gallai tree with 15 blocks.
}
\label{gallai-tree-pic}
\end{figure}
\end{center}

\begin{thm}\label{d0Classify}
A connected graph is degree-choosable if and only if it is not a Gallai tree.
\end{thm}

A minimal counterexample to Brooks' Theorem is regular, and hence
degree-choosable if and only if $\Delta$-choosable.  So
Brooks' Theorem follows immediately from Theorem \ref{d0Classify},
since the only $\Delta$-regular Gallai tree is $K_{\Delta+1}$.
We give two proofs of Theorem \ref{d0Classify}.
The first uses a structural lemma from \erdos, Rubin, and Taylor
\cite{ErdosRT79} known as ``Rubin's Block Lemma,'' but the second does not.

\begin{lem}
No Gallai tree is degree-choosable.
\label{gallai-tree}
\end{lem}
\begin{proof}
Assign disjoint lists to the blocks as follows.
For each block $B$, let $\LB$ be a list of size 2 if $B$ is an odd cycle and
size $k$ if $B$ is $K_{k+1}$.  The list for each vertex $v$
is the union of the lists for blocks containing it:
$L(v)=\cup_{B\ni v}\LB$.  Note that $|L(v)|=d(v)$ for all $v$.

We show that $G$ has no coloring from these lists, by induction on the
number of blocks.  For the base case, $G$ is complete or an odd cycle, and
the proposition clearly holds.  Otherwise, let $B$ be an endblock of $G$, and
let $v\in B$ be a cutvertex.  If $G$ has an $L$-coloring, then each vertex in
$B-v$ gets colored from $\LB$, so each color in $\LB$ appears on a neighbor of
$v$.  Now $G\setminus (B-v)$ is again a Gallai
tree, with lists as specified above; by hypothesis it has no good coloring from
its lists.
\end{proof}


The following lemma has many different proofs 
\cite{ErdosRT79, Entringer1985367, HladkyKS10, raberndiss}.  We follow the
presentation of Hladk\'{y}, Kr\'{a}\soft{l}, and Schauz \cite{HladkyKS10}.
Although it is known as Rubin's Block~\cite{ErdosRT79}, the result was implicit
in the much earlier work of Gallai~\cite{Gallai63a,Gallai63b} and Dirac.

\begin{Rubin} 
\label{rubin}
If $G$ is a $2$-connected graph that is not complete and not an odd cycle,
then $G$ contains an induced even cycle with at most one chord.
\end{Rubin}
\begin{proof}
Let $G$ be a 2-connected graph that is neither complete nor an odd cycle. 
Since $G$ is not complete, it has a minimal cutset $S$, with $|S|\ge 2$.  
Choose $u,v\in S$ and let $C$ be a cycle formed from the union of shortest
paths $P_1$ and $P_2$ joining $u$ and $v$ in two components of $G\setminus S$
(see Figure \ref{block-lemma}). 
Now $C$ has at most one chord, the edge $uv$.  If $C$ has even length, then we
are done.  If $C$ has odd length, then one of the paths joining $u$ and $v$ in
$C$ has odd length; call it $P$.  If $uv$ is present, then $P+uv$ is a
chordless even cycle.  Thus, $uv$ is absent and $C$ is an induced odd cycle
of length at least 5.
Since $G$ is not an odd cycle, there exists $w\in V(G\setminus C)$.

Suppose first that no vertex $w\in V(G\setminus C)$ has two neighbors on $C$.
Since $G$ is 2-connected, there is a shortest path $R$ with endpoints on $C$
and interior disjoint from $C$, and $R$ has length at least 3.  
Now $V(C)\cup V(R)$ induces
two 3-vertices with three vertex disjoint paths between them.  Two 
paths have the same parity, so together they induce a chordless even cycle.

So instead some vertex $w\in V(G\setminus C)$ has two or more neighbors on $C$;
call them $v_1,\ldots, v_k$ (see Figure~\ref{block-lemma}).  The $v_i$ split
$C$ into paths $P_i$ with each $v_i$ the endpoint of two paths.  If any $P_i$
has even length, then $V(P_i) \cup \{w\}$ induces a chordless even cycle.  So
each $P_i$ has odd length; since $C$ has odd length, $k\ge 3$.  If $k>3$, then
$V(P_1)\cup V(P_2)\cup\{w\}$ induces an even cycle with one chord.  If $k=3$,
then some path $P_i$, say $P_3$ by symmetry, has length at least 3, since $C$
has length at least 5.  Now again $V(P_1)\cup V(P_2) \cup\{w\}$ induces an even
cycle with exactly one chord.
\end{proof}

\begin{center}
\begin{figure}[h!tb]
\begin{tikzpicture}[line width=.03cm, scale=.7]
\tikzstyle{vert}=[draw,shape=circle,fill=black,minimum size=4pt, inner sep=0pt] 
\tikzstyle{small}=[draw,shape=circle,fill=black,minimum size=2pt, inner sep=0pt] 

\node[vert,anchor=east,label=below:{$u$}] at (0,0) (u) {};
\node[vert,anchor=east,label=above:{$v$}] at (0,3) (v) {};

\draw (v) .. controls ([xshift=-3cm,yshift=1cm] v)
      and ([xshift=-3cm, yshift=-1cm] u) .. (u);
\draw (u) .. controls ([xshift=3cm,yshift=-1cm] u)
      and ([xshift=3cm, yshift=1cm] v) .. (v);

\draw (u) to [bend left=20] (v);
\draw (v) to [bend left=20] (u);
\draw[dashed] (v) to  (u);

\draw[line width=.07cm] (u) to [bend left=60] node[auto] {$P_1$} (v);
\draw[line width=.07cm] (v) to [bend left=60] node[auto] {$P_2$} (u);

\begin{scope}[xshift=5.5cm, yshift=.25cm, line width =.03cm, scale=.7]
\draw (0,0) node[vert,label=left:{$v_1$}] (v1) {}
   -- (1,0) node[small] {}
   -- (2,0) node[small] {}
   -- (3,0) node[vert] (v3) {}
   -- (4,0) node[small] {}
   -- (5,0) node[small] {}
   -- (6,0) node[vert] (v6) {}
   -- (7,0) node[small] {}
   -- (8,0) node[small] {}
   -- (9,0) node[vert,label=right:{$v_k$}] (vk) {};
\draw (v1) -- (4.5,3) node[vert,label=above:{$w$}] (w) {} -- (vk);
\draw (v3) -- (w) -- (v6); 
\draw (v1) .. controls ([yshift=-1.5cm] v1)
      and ([yshift=-1.5cm] vk) .. (vk);


\end{scope}
\end{tikzpicture}

\caption{The figure on the left shows an induced cycle $C$ formed from $P_1$
and $P_2$.  The figure on the right shows $w$, not on $C$, and its neighbors on
$C$.
}
\label{block-lemma}
\end{figure}
\end{center}

The following proof is similar to Lovasz's proof of Brooks' Theorem, which we
saw in Section~\ref{lovasz}.  In the list coloring context, two parts of that
proof break down: (i) nonadjacent neighbors of a common vertex $v$ need not
have a common color and (ii) we cannot permute colorings on different blocks to
agree on a cutvertex.  The first problem has an easy solution, but the
second is more serious.  If in our induced path $uvw$ either $u$ or $w$ is a
cutvertex, then when we color $u$ and $w$ first we disconnect the graph of
uncolored vertices.  So rather than coloring toward a subgraph that we colored
first, we instead color toward a subgraph that we can color last---any
degree-choosable subgraph will do.

\begin{proof}[Proof 1 of Theorem \ref{d0Classify}]
Lemma~\ref{gallai-tree} shows that no Gallai tree is degree-choosable.
So now suppose $G$ is not a Gallai tree.  By \hyperref[rubin]{Rubin's Block
Lemma}, $G$ has an induced even cycle with at most one chord; call this
subgraph $H$.  We can greedily color the vertices of $G-H$ in order of
decreasing distance from $H$, since each vertex has an uncolored neighbor when
we color it.  Now we extend the coloring to $H$.  If $H$ is an even cycle, we
use Lemma ~\ref{cycle-lemma}.  So assume instead that $H$ is an even cycle with
one chord; label the vertices $v_1,\ldots, v_n$ around the cycle so that
$d_H(v_1)=3$.  Since $v_1$ has 3 colors, we color it with some color not
available for $v_n$.  Now we greedily color the vertices in order of increasing
index.
\end{proof}

Many proofs of list coloring theorems can be easily extended to prove their
analogues for online list coloring.  
The classification of degree-choosable graphs illustrates this idea well. 
The analogue of degree-choosable is
\emph{degree-paintable}, and a connected graph is degree-paintable precisely
when it is not a Gallai tree.  Suppose that $H$ is a connected degree-paintable
induced subgraph of $G$.  Let $\sigma$ be an order of the vertices by
decreasing distance from $H$.  If $S_k$ denotes the vertices available on round
$k$, then the algorithm greedily forms a maximal independent set $I_k$ by
adding vertices from $(V(G)-V(H))\cap S_k$ in the order $\sigma$.  For the
vertices in $S_k$ with no neighbor in $I_k$, the algorithm then plays on $H$ the
strategy that shows it is degree-paintable.  This produces a valid coloring.
To complete the classification of degree-paintable graphs, we need only verify
that every even cycle with at most one chord is degree-paintable.
\smallskip

Alon and Tarsi \cite{AlonT92} developed a powerful tool, which gives an
alternate short proof of Theorem \ref{d0Classify} from 
\hyperref[rubin]{Rubin's Block Lemma}.  
A
subgraph $H$ of a directed graph $D$ is \emph{Eulerian} if in $H$ each vertex
$v$ has indegree $d_H^-(v)$ equal to outdegree $d_H^+(v)$.  Let $EE$ and $EO$ denote the
sets of Eulerian subgraphs of $D$ where the number of edges is even and odd,
respectively.

\begin{AT-thm}
Let $D$ be an orientation of a graph $G$, and let $L$ be a list assignment
such that $|L(v)|\ge 1 + d^+(v)$ for all $v$.
If $|EE|\ne |EO|$, then $G$ is $L$-colorable.
\end{AT-thm}

Let $G$ be a graph that is not a Gallai Tree, and let $H$ be an induced even cycle
with at most one chord in $G$, as guaranteed by 
\hyperref[rubin]{Rubin's Block Lemma}.  
Order the vertices outside of $H$ by increasing
distance from $H$ (breaking ties arbitrarily); call this order $\sigma$.  
Now orient each edge $uv$ as $u\to v$ if $u$ precedes $v$ in $\sigma$.
Orient the cycle edges in $H$ consistently, and
orient the chord, if it exists, arbitrarily.  It is easy to see that each
vertex $v$ has indegree as least 1, and hence outdegree at most $d(v)-1$.  So
the Alon--Tarsi Theorem implies that $G$ is degree-choosable if $|EE|\ne |EO|$.
Every Eulerian subgraph $J$ must be a subgraph of $H$, for otherwise the vertex
of $J$ that comes last in $\sigma$ has outdegree 0.  If $H$ is an even cycle,
then $|EO|=0$ and $|EE|=2$, since both the empty subgraph and all of $H$ are in
$EE$.  If $H$ is a cycle with a chord, then either $|EE|=3$ and $|EO|=0$ or
else $|EE|=2$ and $|EO|=1$.  In both cases, the Alon--Tarsi Theorem shows that
$G$ is degree-choosable.  This is essentially the proof given by
Hladk\'{y}, Kr\'{a}l, and Schauz~\cite{HladkyKS10}.

After our first proof of Theorem~\ref{d0Classify}, we outlined how to
extend the result to characterize degree-paintable graphs.  However, we
omitted the tedious proof that an even cycle with at most one chord is
degree-paintable. 
One advantage of the Alon--Tarsi proof of Theorem~\ref{d0Classify} 
is that it extends easily to paintability.  
Schauz proved an analogue of the Alon--Tarsi Theorem for {online}
list coloring\footnote{Let $D$ be an orientation of a graph $G$, and let $f$ be
a list size assignment such that $f(v)\ge 1 + d^+(v)$ for all $v$.
If $|EE|\ne |EO|$, then $G$ is online $f$-choosable.}, so in fact the proof of
the paintability result is nearly identical.

We extend this idea \cite{CranstonR13+paint} to show that $G^2$ is online
$(\Delta^2-1)$-choosable unless $\omega(G^2)\ge \Delta^2$; here $G^2$ is
formed from $G$ by adding edge $uv$ for each pair $u$ and $v$ at distance 2 in
$G$.  Our approach in that proof 
is quite similar to method used above to prove Theorem~\ref{d0Classify} via the
Alon--Tarsi Theorem.  Applying these techniques to $G^2$
gives $d^-(v)\ge 2$ for all $v$.  Now however, we seek a subgraph $H$
that is online ``degree-1''-choosable, that is, $H$ is online $f$-choosable,
where $f(v)=d(v)-1$.  So the bulk of the work lies in showing that every square
graph contains either such an induced subgraph or else a large clique.
\bigskip

We now conclude our digression into online list coloring and the Alon--Tarsi
Theorem, and we finish the section with a second proof of
Theorem~\ref{d0Classify}.
%
Kostochka, Stiebitz, and Wirth~\cite{BrooksExtended} gave a short proof of
Theorem \ref{d0Classify}, which we reproduce below.
Further, they extended the result to hypergraphs.  

\begin{proof}[Proof 2 of Theorem \ref{d0Classify}]
Lemma~\ref{gallai-tree} shows that no Gallai tree is degree-choosable.
So now suppose there exists a graph that is not a Gallai tree, but is also not
degree-choosable.  Let $G$ be a minimum such graph.
Since $G$ is not degree-choosable, no induced subgraph $H$ of $G$ is
degree-choosable (if such an $H$ exists, then we color $G-H$ greedily towards $H$, and
extend the coloring to $H$ since it is degree-choosable). 
Hence for any $v \in V(G)$ that is not a cutvertex, $G-v$ must be a Gallai
tree by minimality of $|G|$.  

If $G$ has more than one block, then for endblocks $B_1$ and $B_2$, choose
noncutvertices $w\in B_1$ and $x\in B_2$.  By the minimality of $|G|$, both
$G-w$ and $G-x$ are Gallai trees, so every block of $G$ is either complete or
an odd cycle, and thus $G$ is a Gallai tree.  So instead $G$ has only one block,
that is, $G$ is $2$-connected.
Further, $G-v$ is a Gallai tree for all $v \in V(G)$.

Now let $L$ be a list assignment on $G$ such that $|L(v)| = d(v)$ for all $v
\in V(G)$ and $G$ is not $L$-colorable.  Suppose two vertices in $G$ get
different lists.  Since $G$ is connected, we have adjacent $v,w \in V(G)$
such that $L(v) - L(w) \ne \emptyset$.   Pick $c \in L(v) - L(w)$ and color $v$
with $c$.  Now we can finish by coloring in order of decreasing distance from
$w$ in $G-v$.  So instead $L(v) = L(w)$ for all $v,w \in V(G)$; in particular,
$G$ is regular.

Pick $v \in V(G)$ and consider the Gallai tree $G-v$.  Since $G$ is regular and
$2$-connected, $v$ must be adjacent to all noncutvertices in all endblocks of $G-v$.  So, if $G-v$ has at least two endblocks, then
$d(v) \ge 2(\Delta(G)-1)$. Since $ 2(\Delta(G)-1) > \Delta(G)$ when $\Delta(G)
\ge 3$, we must have $\Delta(G) = 2$. Since $G$ is not $L$-colorable it is not
$2$-colorable and hence is an odd cycle, a contradiction.  Therefore
$G-v$ has only one endblock and thus $G$ is complete, again a
contradiction.  
\end{proof}

\section{Further Directions}
In this final section, we conclude our survey by discussing two conjecture that
strengthen Brooks' Theorem.
Determining the chromatic number of a graph is well-known to be NP-hard.
However, Brooks' Theorem shows that determining whether a graph has
$\chi=\Delta+1$ is easy.  If $\Delta=2$, look for odd cycles; otherwise, look
for a $K_{\Delta+1}$.
It is natural to ask how close $t$ must be to $\Delta$ so that we can easily 
check whether $\chi=t$.  Emden-Weinert, Hougardy, and Kreuter
\cite{Emden-WeinertHK98} gave the following lower bound on this threshold.

\begin{thm}
For any fixed $\Delta$, deciding whether a graph $G$ 
with maximum degree $\Delta$ 
has a $(\Delta+1-k)$-coloring is NP-complete for any $k$ such that
$k^2+k>\Delta$, when $\Delta+1-k\ge 3$. 
\end{thm}

Molloy and Reed \cite{MolloyR01STOC} then proved a matching upper bound, for
sufficiently large $\Delta$.

\begin{thm}
For any fixed sufficiently large $\Delta$, deciding whether a graph $G$ 
with maximum degree $\Delta$ 
has a $(\Delta+1-k)$-coloring is in $P$ for every $k$ such that
$k^2+k\le \Delta$.
\end{thm}

Further, they conjectured that the same result holds for all values of $\Delta$.
(Section 15.4 of Molloy and Reed \cite{MolloyR-GCPM} has more on this question.)
Viewed in this framework, Brooks' Theorem describes the case $k=1$.
Now we consider the case $k=2$.  In 1977, Borodin and
Kostochka~\cite{BorodinK77} posed the following conjecture.

\begin{BK}
If $G$ has $\Delta\ge 9$ and $\omega \le \Delta-1$, then
$\chi\le \Delta-1$.
\end{BK}
If true, the conjecture is best possible in two ways.  
First, even when we require $\omega\le \Delta-2$, we cannot conclude 
$\chi\le \Delta-2$.  For example, form $G$ from a disjoint 5-cycle and 
$K_{\Delta-4}$ by adding all edges with one endpoint in each graph (this is the
\emph{join} of $C_5$ and $K_{\Delta-4}$).  Now $\omega=\Delta-2$, but
$\chi=\Delta-1$, since every proper coloring uses at least 3 colors on the
5-cycle and cannot reuse any color on the clique.

\begin{figure}[!ht]
\begin{center}
\begin{tikzpicture}[rotate=18,scale=.4,minimum size = 2pt, inner sep=0pt]
   \GraphInit[vstyle=Hasse]

\tikzset{VertexStyle/.style = {shape = circle,fill = black,minimum size = 5pt,inner sep=0pt}}

\grCycle[prefix=a,RA=2.75]{5}
\grCycle[prefix=b,RA=4.25]{5}
\grCycle[prefix=c,RA=5.75]{5}
\EdgeMod{a}{b}{5}{1}
\EdgeMod{b}{c}{5}{1}
\EdgeMod{a}{c}{5}{1}
\EdgeMod{a}{b}{5}{4}
\EdgeMod{b}{c}{5}{4}
\EdgeMod{a}{c}{5}{4}
\EdgeIdentity{a}{b}{5}
\EdgeIdentity{b}{c}{5}
\Edge  [style ={out = 215, in=-20}](c0)(a0)
\Edge  [style ={out = 287, in=52}](c1)(a1)
\Edge  [style ={out = 359, in=124}](c2)(a2)
\Edge  [style ={out = 71, in=196}](c3)(a3)
\Edge  [style ={out = 143, in=268}](c4)(a4)
\end{tikzpicture}
\end{center}
\caption{A construction showing that the hypothesis $\Delta\ge 9$ in the
Borodin--Kostochka Conjecture is necessary and best possible.}
\label{BK-pic}
\end{figure}

Second, the lower bound on $\Delta$ cannot be reduced, as shown by the
following construction.  Form $G$ from five disjoint copies of $K_3$, say
$D_1,\ldots,D_5$, by adding edge $uv$ if $u\in D_i$, $v\in D_j$, and
$|i-j|\equiv1\bmod 5$.  This graph is 8-regular with $\omega=6$.
Each color class has size at most 2, 
so $\chi=\ceil{15/2}=8$. 
Thus 
$\chi=\Delta$, but $\omega=\Delta-2$.
For $\Delta\le 8$, various other examples are known~\cite{CranstonR13claw} 
where $\chi=\Delta$ and $\omega<\Delta$.

Reed \cite{Reed99brooks} proved the Borodin--Kostochka Conjecture when
$\Delta\ge 10^6$ and the present authors \cite{CranstonR13claw} proved it for
claw-free graphs (those where no vertex has three pairwise nonadjacent
neighbors).  Although this question remains open, stronger versions of the
conjecture are believed true. 
Already in 1977, Borodin and Kostochka were convinced that the same upper bound
holds for the list chromatic number.  Recently, we conjectured
\cite{CranstonR13+paint} that the bound holds even for the online list chromatic
number.


Reed posed the following conjecture, which is along similar lines, but much
more far-reaching.

\begin{Reed}
Every graph $G$ satisfies $\chi \le \ceil{\frac{\omega+\Delta+1}2}$.
\end{Reed}

In 1998, Reed proved \cite{Reed98omega} that
there exists a positive $\epsilon$ such that 
$\chi \le \ceil{\omega\epsilon + (\Delta+1)(1-\epsilon)}$. 
The original conjecture is
that this upper bound holds when $\epsilon = \frac12$.  In 2012, King and
Reed~\cite{KingR12ep} gave a much shorter proof of the same result.  A key
ingredient of their proof is the result of King \cite{KingHitting} that every
graph with $\omega > \frac{2(\Delta+1)}3$ has a hitting set (recall that a
hitting set is independent and intersects every maximum clique).  
About the same time, they used the Claw-free Structure Theorem of Chudnovsky
and Seymour to prove that Reed's Conjecture holds for all claw-free
graphs~\cite{KingR12claw-free}.
Section 21.3 of Molloy and Reed~\cite{MolloyR-GCPM} gives further evidence for
Reed's Conjecture by showing that the desired upper bound holds for the
fractional chromatic number, even without rounding up.

We conclude this section by showing that Reed's Conjecture is best
possible, using random graphs.  Specifically, we show that if $\epsilon >
\frac12$, then the bound $\chi \le \ceil{\omega\epsilon +
(\Delta+1)(1-\epsilon)}$ fails for some graph $G$.  The proof we present is
from the end of~\cite{Reed98omega}.
Kostochka~\cite{Kostochka84} also showed this using the explicit examples that
Catlin \cite{catlin1979hajos} constructed to disprove the Haj{\'o}s Conjecture. 
Let $H_t
= t \cdot C_5$ (i.e.~$C_5$ where each edge has multiplicity $t$) and let $G_t$ 
be the line graph of $H_t$; Figure~\ref{BK-pic} shows $G_3$.
Catlin showed that for odd $t$ we have $\chi(G_t) = \frac{5t +
1}{2}$, $\Delta(G_t) = 3t - 1$, and $\omega(G_t) = 2t$.  So, for any $\epsilon
> \frac12$, we can choose $t$ large enough to make the bound fail.

We will randomly construct a graph $H$ on $n$
vertices such that $\chi(H)\ge \frac12n-n^{3/4}$ and $\omega(H)\le 8n^{3/4}\log
n$.  When we form $G$ by joining $H$ to $K_{\Delta+1-n}$, we see that the
desired bound fails for $G$ when $\epsilon > \frac12$.
Let $H$ be a random graph on $n$ vertices, where each edge appears
independently with probability $p$ and let $p=1-n^{-3/4}$.  The expected number
of cliques of size $k$ is ${n \choose k}p^{k\choose 2}$.  So when $k\ge
8n^{3/4}\lg n$, the expected number of $k$-cliques is arbitrarily
small for sufficiently large $n$.
\begin{align*}
{n \choose k}p^{k\choose 2} & = 
{n \choose k}\left(1-n^{-3/4}\right)^{k\choose 2} \\
& \le n^k \left(1-n^{-3/4}\right)^{k^2/4} \\
& \le 2^{k\lg n}\left(e^{-n^{-3/4}}\right)^{k^2/4}\\
& \le e^{k\lg n-k^2 n^{-3/4}/4}\\
& \le e^{-k \lg n}.
\end{align*}
The expected number of independent sets of size 3 is
${n\choose 3}p^3 = {n\choose 3}(1-(1-n^{-3/4}))^3\le\frac16n^{3/4}$.
By deleting one vertex from each independent set of size 3, we get a graph $H$
on $n-\frac16n^{3/4}$ vertices with independence number 2.  So $\chi(H)\ge
\frac{n}2-\frac1{12}n^{3/4}$ and $\omega(H)\le 8n^{3/4}\lg n$.
Now $\chi(G) \ge (\Delta+1-n) +
(\frac12n-\frac{1}{12}n^{3/4})=\Delta+1-\frac{1}{2}n-\frac{1}{12}n^{3/4}$.  
Similarly $\omega(G)\le (\Delta+1-n)+8n^{3/4}\lg n$.

\section*{Acknowledgments}
We thank many people for helpful comments on earlier versions of this paper.
These include
Katie Edwards,
Richard Hammack,
Bobby Jaeger,
Andrew King,
Sasha Kostochka,
Bernard Lidicky,
Bojan Mohar,
Riste \v{S}krekovski,
Carsten Thomassen,
Helge Tverberg,
and
Gexin Yu.
In particular, Bobby Jaeger, Katie Edwards, and Bernard Lidicky
gave very detailed constructive feedback.
We also should mention a recent preprint~\cite{StiebitzT13+} by Michael
Stiebitz and Bjarne Toft entitled ``Brooks' Theorem''; their survey has some
similarities to this one, but a different tone and focus.

\nocite{ponstein1969new}
\nocite{dirac1957theorem}
\bibliographystyle{abbrvplain}
\bibliography{GraphColoring}
\end{document}

%% file: cycle-figs.tex
%
\begin{center}
\begin{figure}[ht]
\begin{tikzpicture}[scale = 1]
\tikzstyle{VertexStyle}=[shape = circle, minimum size = 22pt, 
inner sep = 1.2pt, draw]

\Vertex[x = 0, y = 0, L = \tiny {$n$}]{v1}
\Vertex[x = 1, y = 0, L = \tiny {$1$}]{v2}
\Vertex[x = -1, y = 0, L = \tiny {$n-1$}]{vn}
\Vertex[x = 1, y = 3, L = \tiny {$k-1$}]{vk-1}
\Vertex[x = 0, y = 3, L = \tiny {$k$}]{vk}
\Vertex[x = -1, y = 3, L = \tiny {$k+1$}]{vk+1}
\Edge [style = {out=0, in=0}](v2)(vk-1)
\Edge (vk-1)(vk)
\Edge (vk)(v1)
\Edge (v1)(vn)
\Edge [style = {out=180, in=180}](vn)(vk+1)
\Edge [style = {dashed}](v2)(vk+1)

\begin{scope}[xshift=2.0in]
\Vertex[x = 0, y = 0, L = \tiny {$n$}]{v1}
\Vertex[x = 1, y = 0, L = \tiny {$1$}]{v2}
\Vertex[x = 2, y = 0, L = \tiny {$2$}]{v3}
\Vertex[x = -1, y = 0, L = \tiny {$n-1$}]{vn}
\Vertex[x = 2, y = 3, L = \tiny {$k-2$}]{vk-2}
\Vertex[x = 1, y = 3, L = \tiny {$k-1$}]{vk-1}
\Vertex[x = 0, y = 3, L = \tiny {$k$}]{vk}
\Vertex[x = -1, y = 3, L = \tiny {$k+1$}]{vk+1}
\Edge (v1)(vk)
\Edge (vk)(vk-1)
\Edge [style = {out=250,in=50}](vk-1)(vn)
\Edge [style = {out=180,in=180}](vn)(vk+1)
\Edge [style = {out=310,in=110}](vk+1)(v2)
\Edge (v2)(v3)
\Edge [style = {out=0,in=0}](v3)(vk-2)
\Edge [style = {dashed}](vk-2)(v1)

\begin{scope}[xshift=2.4in]
\Vertex[x = 0, y = 0, L = \tiny {$n$}]{v1}
\Vertex[x = 1, y = 0, L = \tiny {$1$}]{v2}
\Vertex[x = -1, y = 0, L = \tiny {$n-1$}]{vn}
\Vertex[x = 1, y = 3, L = \tiny {$k-1$}]{vk-1}
\Vertex[x = 0, y = 3, L = \tiny {$k$}]{vk}
\Vertex[x = -1, y = 3, L = \tiny {$k+1$}]{vk+1}
\Edge [style = {out=0, in=0}](v2)(vk-1)
\Edge (vk-1)(vk)
\Edge (vk)(v1)
\Edge [style = {out=130, in=50}](v2)(vn)
\Edge [style = {out=180, in=180}](vn)(vk+1)
\Edge [style = {dashed}](v1)(vk+1)

%
\end{scope}
\end{scope}
\end{tikzpicture}
\caption{Key steps in the proof that (4) $\implies$ (1): On the left, $G$
must contain $v_1v_{k+1}$.  In the center, $G$ must contain $v_{n}v_{k-2}$.
A symmetric argument shows that $G$ must contain $v_{n}v_{k+2}$.
On the right, $G$ must contain $v_{n}v_{k+1}$.  
}
\label{cycle-figs}
\end{figure}
\end{center}